\RequirePackage{fix-cm}
\documentclass[compress,preprint]{elsarticle}
\biboptions{sort&compress}

\usepackage[abs]{overpic}
\usepackage{graphicx}
\usepackage[fontsize=9pt]{fontsize}
\usepackage{amsfonts, color, amssymb}
\usepackage{amsmath}
\usepackage{amsthm}
\usepackage{epsfig}
\usepackage{psfrag}
\usepackage{bm}
\usepackage{paralist}
\usepackage{algorithmicx}
\usepackage{algpseudocode}
\usepackage{appendix}
\usepackage{caption}
\usepackage{hyperref}
\usepackage{xurl}
\hypersetup{hidelinks}
\usepackage{color}
\usepackage{booktabs}
\usepackage{tabularx}
\usepackage{subcaption}
\usepackage{siunitx}
\usepackage{microtype}
\usepackage{needspace}
\usepackage{placeins}
\usepackage{xcolor}
\usepackage{tikz}
\usetikzlibrary{arrows.meta}


\usepackage[misc]{ifsym}

\newcommand{\dt}{\Delta t}
\newcommand{\gradop}{\operatorname{grad}}
\newcommand{\curlop}{\operatorname{curl}}
\newcommand{\divop}{\operatorname{div}}
\newcommand{\curl}{\nabla\times}
\newcommand{\curlh}{\nabla_h\times}
\newcommand{\diver}{\nabla\cdot}
\newcommand{\norm}[1]{\left\|#1\right\|}
\newcommand{\Eh}{E_h}
\newcommand{\Fh}{F_h}
\newcommand{\Vh}{V_h}
\newcommand{\Vper}{V_{\mathrm{per}}}
\newcommand{\Eper}{E_{\mathrm{per}}}
\newcommand{\Fper}{F_{\mathrm{per}}}
\newcommand{\Dp}{D_+}
\newcommand{\Dm}{D_-}

\newcommand{\Qh}{Q_h}

\usepackage[a4paper, total={6.55in,9.15in}]{geometry}

\makeatletter
\patchcmd{\ps@pprintTitle}
  {\let\@evenfoot\@oddfoot}
  {\apptocmd{\@oddfoot}{\unskip}{}{}\let\@evenfoot\@oddfoot}
  {}{}
\makeatother

\allowdisplaybreaks[4]
\newtheorem{theorem}{Theorem}[section] %
\newtheorem{remark}{Remark}[section]
\newcounter{algorithm}[section]
\renewcommand{\thealgorithm}{\thesection.\arabic{algorithm}}
\newenvironment{algorithm}[1][]{%
  \par\addvspace{8pt plus 2pt minus 2pt}%
  \Needspace{5\baselineskip}%
  \refstepcounter{algorithm}%
  \normalfont\setlength{\parindent}{0pt}\setlength{\parskip}{0pt}%
  \postdisplaypenalty=10000\relax
  {\hrule height 0.8pt\relax}%
  \nobreak\vskip 3pt\relax
  \noindent\textbf{Algorithm~\thealgorithm}%
  \if\relax\detokenize{#1}\relax\else\quad #1.\fi
  \par\nobreak\vskip 3pt\relax
  {\hrule height 0.4pt\relax}%
  \nobreak\vskip 5pt\relax
  \ignorespaces
}{%
  \par\nobreak\vskip 5pt\relax
  {\hrule height 0.4pt\relax}%
  \par\addvspace{8pt plus 2pt minus 2pt}%
  \ignorespacesafterend
}
\newtheorem{proposition}[theorem]{Proposition} %

\newtheorem{corollary}[theorem]{Corollary}
\AtBeginEnvironment{proof}{%
\setlength{\abovedisplayskip}{5pt plus 1pt minus 1pt}
\setlength{\belowdisplayskip}{5pt plus 1pt minus 1pt}
\setlength{\abovedisplayshortskip}{3pt plus 1pt minus 1pt}
\setlength{\belowdisplayshortskip}{3pt plus 1pt minus 1pt}
\setlength{\jot}{2pt}
}

\def\3bar{{|\hspace{-.02in}|\hspace{-.02in}|}}
\allowdisplaybreaks[4]

\allowdisplaybreaks 

\numberwithin{equation}{section}

\begin{document}

\begin{frontmatter}

\title{A linear dual-field finite element method for  incompressible MHD, part I: cross-helicity and Alfv\'enic alignment in periodic domains}

\author[1]{Shusheng Li}
\ead{ssli22@mails.jlu.edu.cn}

\author[1]{Siyuan Tong\corref{cor1}}
\ead{tongsy22@mails.jlu.edu.cn}
\cortext[cor1]{Corresponding author}

\author[1]{Qilong Zhai}
\ead{zhaiql@jlu.edu.cn}

\author[1]{Qian Zhang}
\ead{qzhang25@jlu.edu.cn}

%% Author affiliation
\affiliation[1]{organization={School of Mathematics, Jilin University},%Department and Organization
	% addressline={}, 
	city={Changchun},
	postcode={130012}, 
	state={Jilin},
	country={People's Republic of China}}

\begin{abstract}
Cross helicity, together with total energy, determines the Els\"asser energies and thus quantifies departures from Alfv\'enic alignment in incompressible MHD. Building on the dual-field method in [Mao and Xi, \emph{J. Comput. Phys.}, 2025], we construct a linear finite element method for periodic domains that reproduces the cross-helicity balance. We connect the discrete energy and cross-helicity conservation laws to Alfv\'enic alignment through a nonnegative Els\"asser-energy representation, establishing preservation of the corresponding exact relation when either Els\"asser energy is initially zero in the unforced ideal case. We  reassign the two discrete approximations of the magnetic field to the integer and half-integer time levels and redesign the Lorentz and induction couplings within the Mao-Xi formulation. The resulting scheme uses linear staggered updates and  retains the energy balances, strong and weak divergence-free properties, and magnetic- and fluid-helicity balances. Numerical experiments illustrate the role of the cross-helicity balance in preserving an aligned state in the ideal case and in alignment dynamics under weak dissipation.
\end{abstract}


\end{frontmatter}

\section{Introduction}
\label{sec-introduction}

Incompressible magnetohydrodynamics (MHD) couples the evolution of a conducting fluid and a magnetic field through magnetic induction and the Lorentz force \cite{davidson2001}.
The energy balance constrains the evolution of the combined kinetic and magnetic energy, but does not determine the correlation between velocity and magnetic field.
Cross helicity measures the global correlation between these fields.
Velocity--magnetic alignment is important for understanding MHD turbulence and solar-wind fluctuations, particularly because local directional alignment can weaken nonlinear interactions \cite{matthaeus2008}.
We develop a linear dual-field finite element method that reproduces the cross-helicity balance and study the implications of this balance for velocity--magnetic alignment.

For velocity $\bm u$ and magnetic field $\bm B$ on a periodic domain $\Omega$, cross helicity is $H_c=\int_\Omega\bm u\cdot\bm B\,\mathrm{d}x$.
Its connection with alignment can be expressed through the Els\"asser variables $\bm z_\pm=\bm u\pm\sqrt{s}\,\bm B$, where $s>0$ is the magnetic coupling parameter \cite{elsasser1950}.
The corresponding Els\"asser energies $\tfrac12\|\bm z_\pm\|_{L^2(\Omega)}^2$ have a sum proportional to total energy and a difference proportional to cross helicity.
In the unforced ideal problem, simultaneous conservation of energy and cross helicity therefore preserves both Els\"asser energies separately.
If either vanishes initially, the corresponding exact Alfv\'enic relation $\bm u=\pm\sqrt{s}\,\bm B$ is preserved.
In viscous and resistive MHD, the relative decay of energy and cross helicity can drive relaxation towards an Alfv\'enically aligned state for suitable initial data \cite{dobrowolny1980,ting1986}, motivating numerical methods that reproduce both balance laws.

Incompressible MHD also satisfies differential constraints and other helicity balance laws.
The velocity and magnetic field satisfy $\diver\bm u=\diver\bm B=0$, while $\diver\bm J=0$ follows from $\bm J=\curl\bm B$ and the identity $\diver\curl=0$.
Magnetic-divergence errors can introduce spurious magnetic forces \cite{MR1764250}, and the discrete treatment of incompressibility affects pressure robustness \cite{MR3683678}.
Magnetic helicity carries information about magnetic-field linkage and knottedness and constrains magnetic relaxation and minimum-energy states \cite{MR1612569,MR96542,taylor1974}.
For zero-mean periodic fields, conservation of nonzero magnetic helicity gives a positive magnetic-energy lower bound through the Arnold inequality \cite{MR1612569}.
Fluid helicity measures velocity--vorticity correlation and carries information about vortex-line linkage and knottedness \cite{moffatt1969}.
Its balance includes the Lorentz-force contribution in MHD and reduces to conservation in the inviscid, unforced hydrodynamic limit, providing a further constraint on vortical evolution.

Finite element exterior calculus organizes finite element spaces into discrete de Rham complexes and provides commuting projections and discrete Hodge decompositions for retaining the underlying differential identities \cite{MR2269741,MR2594630,MR2009375}. This framework has supported MHD methods with exactly divergence-free magnetic fields \cite{MR3599563}, exactly divergence-free velocity and magnetic approximations \cite{MR3786791}, and formulations that also retain the divergence-free property of the current density \cite{MR4199376}. For magnetic relaxation, He et al.\ \cite{MR5042378} combined finite element de Rham complexes with an auxiliary magnetic field defined by an $L^2$ projection into the $H(\operatorname{curl})$ space of the complex. This construction preserves magnetic helicity and the magnetic-energy dissipation law of the magneto-frictional system. Their discrete Arnold inequality and numerical experiments demonstrate how helicity preservation prevents nonzero-helicity fields from relaxing to a spurious zero state.

Cross-helicity-preserving discretizations are already available. Liu and Wang \cite{liuwang2004} developed finite difference schemes retaining energy and cross helicity at the semidiscrete level for periodic MHD flows with coordinate symmetry. Kraus and Maj \cite{kraus2018} used discrete exterior calculus to construct a finite difference scheme on a staggered Cartesian grid for two-dimensional ideal incompressible MHD that preserves energy and cross helicity. Hu, Lee, and Xu \cite{MR4236017} developed a finite element method retaining the magnetic divergence constraint and the energy, magnetic-helicity, and cross-helicity laws. Beir\~{a}o da Veiga, Hu, and Mascotto \cite{MR4886367} subsequently established error estimates for the associated seven-field spatial semidiscretization. More recently, Yao et al. \cite{yao2026} developed a variable-step, second-order method in Els\"asser variables that preserves energy, cross helicity, and magnetic helicity in the unforced ideal case. Gawlik and Gay-Balmaz \cite{MR4346694} treated variable-density ideal incompressible MHD, preserving energy, magnetic helicity, and the divergence constraints; their cross-helicity conservation result applies to constant density. For Hall MHD, cross helicity is generally no longer conserved even in the ideal limit. Laakmann, Hu, and Farrell \cite{MR4632710} constructed a discretization preserving energy, magnetic helicity, and hybrid helicity in the ideal, unforced limit, where hybrid helicity combines magnetic, cross, and fluid helicities with coefficients determined by the Hall parameter.

For standard MHD, the fully discrete finite element methods in \cite{MR4236017,yao2026,MR4346694} require nonlinear solves. Case et al. \cite{case2010} proposed a linear Scott–Vogelius finite element scheme that retains the energy and cross-helicity balances and yields pointwise divergence-free velocity and magnetic fields.
Our method is motivated by the dual-field approach \cite{palhagerritsma2017}. For periodic incompressible Navier--Stokes flow, Zhang et al.\ \cite{MR4354366} developed a dual-field method that approximates the velocity in $H(\operatorname{curl})$ and $H(\operatorname{div})$ on staggered time grids and preserves mass, energy, and fluid helicity. For incompressible MHD, Mao and Xi \cite{MR4917731} proposed a formulation with two sets of discrete velocity and magnetic fields, one in $H(\operatorname{curl})$ and the other in $H(\operatorname{div})$. Their method enforces the strong and weak divergence-free constraints and satisfies the energy, magnetic-helicity, and fluid-helicity balance laws. Each staggered update requires the solution of a linear mixed system. It does not, however, generally preserve the discrete cross helicity considered here. 

Building on the dual-field formulation of Mao and Xi \cite{MR4917731}, we incorporate the cross-helicity balance while retaining its sequential linear staggered updates, energy balances, and other structural properties. We further investigate how simultaneous conservation of energy and cross helicity constrains Alfv\'enic alignment in this discrete setting.
Two coordinated changes distinguish our method from the Mao--Xi formulation \cite{MR4917731}. First, we change the temporal staggering of the two magnetic field approximations: at integer times, the velocity--magnetic unknown pair changes from $[\bm u_2,\bm B_2]$ to $[\bm u_2,\bm B_1]$, and at half-integer times, from $[\bm u_1,\bm B_1]$ to $[\bm u_1,\bm B_2]$. Second, we evolve the magnetic fields directly rather than their vector potentials, and re-pair the Lorentz and induction couplings accordingly. These changes provide the nonlinear cancellations required by both energy and cross helicity. Each substep remains linear because the magnetic factor and vorticity used as coupling coefficients are already available from the other time grid. We derive two energy balances and a discrete cross-helicity balance. A nonnegative sum-of-squares representation of the discrete Els\"asser energies relates these balance laws to Alfv\'enic alignment: in the unforced ideal case, simultaneous conservation of energy and cross helicity preserves the initial values of these global measures of departure from Alfv\'enic alignment. In particular, if either discrete Els\"asser energy is initially zero, the corresponding exact Alfv\'enic alignment relation is preserved throughout the evolution. Table~\ref{tab:staggered-pairing} details both changes. Subscripts~1 and~2 identify fields in $H(\operatorname{curl})$ and $H(\operatorname{div})$, respectively. In the table, time superscripts are omitted, and bars denote averages over the corresponding substep. The induction factor is the velocity--magnetic cross product entering the magnetic update. 
The method also retains the strong and weak divergence-free properties of velocity, magnetic field, and current density, together with the magnetic- and fluid-helicity balances in their respective physical limits. In the unforced non-resistive case, conserved nonzero magnetic helicity yields positive lower bounds for both staggered energies. 

\begin{table}[!htbp]
\centering
\caption{Staggered fields and couplings.}
\label{tab:staggered-pairing}
\small
\setlength{\tabcolsep}{3pt}
\renewcommand{\arraystretch}{1.18}
\begin{tabular}{@{}lcc ll cc ll@{}}
\toprule
& \multicolumn{4}{c}{Integer step: $t_n\longrightarrow t_{n+1}$}
& \multicolumn{4}{c}{Half-integer step: $t_{n-1/2}\longrightarrow t_{n+1/2}$}\\
\cmidrule(lr){2-5}\cmidrule(l){6-9}
Method
& \shortstack{Unknown\\pair} & \shortstack{Evolved\\variable}
& \shortstack{Lorentz\\force} & \shortstack{Induction\\factor}
& \shortstack{Unknown\\pair} & \shortstack{Evolved\\variable}
& \shortstack{Lorentz\\force} & \shortstack{Induction\\factor}\\
\midrule
Mao--Xi \cite{MR4917731}
& $[\bm u_2,\bm B_2]$ & $\bm A_1$
& $s\overline{\bm J}_1\times\bm B_1$
& $\overline{\bm u}_2\times\bm B_1$
& $[\bm u_1,\bm B_1]$ & $\bm A_2$
& $s\overline{\bm J}_2\times\bm B_2$
& $\overline{\bm u}_1\times\bm B_2$\\
Proposed
& $[\bm u_2,\bm B_1]$ & $\bm B_1$
& $s\overline{\bm J}_2\times\bm B_2$
& $\overline{\bm u}_2\times\bm B_2$
& $[\bm u_1,\bm B_2]$ & $\bm B_2$
& $s\overline{\bm J}_1\times\bm B_1$
& $\overline{\bm u}_1\times\bm B_1$\\
\bottomrule
\end{tabular}
\end{table}

Through numerical comparisons with the Mao--Xi method, we further investigate the role of cross-helicity preservation in Alfv\'enic alignment. In an ideal test with exactly aligned initial data, both methods conserve energy. The proposed method additionally conserves cross helicity, thereby preserving the initial Alfv\'enic alignment, in agreement with the continuous solution. The Mao--Xi method instead departs from its initial state. For the weakly dissipative initial data considered, the proposed method satisfies the cumulative cross-helicity balance and exhibits increasing global and local alignment. This tendency towards Alfv\'enic alignment is qualitatively consistent with the dynamic-alignment regime reported in the two-dimensional study of Ting et al.~\cite{ting1986}. Convergence and conservation tests further assess accuracy and other conservation properties, while an Orszag--Tang experiment examines the evolution of spatial structures. 

The remainder of the paper is organized as follows. Section~\ref{sec:continuous-formulation} introduces the periodic function spaces and the dual-field weak formulation and establishes the continuous divergence constraints and balance laws. Section~\ref{sec:discrete-scheme} constructs the staggered finite element scheme on a periodic discrete de Rham complex and proves its divergence-free properties and discrete energy and helicity balances. It also derives conditions for preserving initially exact Alfv\'enic alignment and energy lower bounds from magnetic helicity. Section~\ref{sec:numerical-experiments} presents convergence tests, comparisons with the Mao--Xi method to examine the role of cross helicity in alignment, further conservation tests, and Orszag--Tang simulations. Section~\ref{sec:conclusions} concludes the paper.

\section{A dual-field mixed weak formulation for the MHD systems}
\label{sec:continuous-formulation}
\subsection{MHD systems and periodic Sobolev spaces}
\label{subsec:continuous-spaces}
Let $\Omega$ be a three-dimensional periodic box, identified with a flat torus, and let $T>0$ denote the final time. The MHD system considered here, written in rotational form, reads
\[
\left\{
\begin{array}{ll}
	\bm{u}_t+\boldsymbol{\omega}\times \bm{u}+\nabla P
	+Re^{-1}\curl\boldsymbol{\omega}
	-s\bm{J}\times \bm{B}
	=\bm{f},
	& \text{in }(0,T)\times\Omega,
	\\[2mm]
	\boldsymbol{\omega}=\curl \bm{u},
	& \text{in }(0,T)\times\Omega,
	\\[2mm]
	\diver \bm{u}=0,
	& \text{in }(0,T)\times\Omega,
	\\[2mm]
	\bm{B}_t+Rm^{-1}\curl \bm{J}
	-\curl(\bm{u}\times \bm{B})
	=\bm{0},
	& \text{in }(0,T)\times\Omega,
	\\[2mm]
	\bm{J}=\curl \bm{B},
	& \text{in }(0,T)\times\Omega,
	\\[2mm]
	\diver \bm{B}=0,
	& \text{in }(0,T)\times\Omega,
	\\[2mm]
	\bm{u}(0)=\bm{u}_0,\qquad
	\bm{B}(0)=\bm{B}_0,
	& \text{in }\Omega.
\end{array}
\right.
\]
Here $\bm u$ is the velocity, $P$ is the Bernoulli pressure, $\boldsymbol\omega=\curl\bm u$ is the vorticity, $\bm B$ is the magnetic field, and $\bm J=\curl\bm B$ is the current density. The body force $\bm f$ and the initial fields $\bm u_0,\bm B_0$ are prescribed. The
numbers $Re$ and $Rm$ are the Reynolds and magnetic Reynolds numbers; hence $Re^{-1}$ and $Rm^{-1}$ are the viscous and resistive coefficients. The parameter $s>0$ weights the magnetic energy and the Lorentz force. No external magnetic forcing is considered. The first equation is the momentum equation. In its rotational form, $\boldsymbol\omega\times\bm u$ is the convection term, $Re^{-1}\curl\boldsymbol\omega$ is the viscous diffusion term, and $s\bm J\times\bm B$ is the Lorentz force. The fourth equation is the magnetic induction equation, in which $Rm^{-1}\curl\bm J$ is the resistive diffusion term and $\curl(\bm u\times\bm B)$ is the induction term.
Unlike the magnetic-potential reformulation in \cite{MR4917731}, we construct the magnetic updates directly from the induction equation for $\bm B$. Magnetic potentials are not time-marching unknowns; they are reconstructed, with a fixed gauge, for the magnetic-helicity analysis. 

We use a periodic de Rham complex to define the function spaces for the MHD formulation. The general differential-form and Hilbert-complex framework is reviewed in \cite[Sec.~2]{MR2269741} and \cite[Secs.~3 and~4]{MR2594630}; the zero-mean reduction below is specific to our periodic setting. Throughout, $(\cdot,\cdot)$ and $\|\cdot\|$
denote the $L^2(\Omega)$ inner product and norm, respectively.
The underlying Sobolev spaces are
\begin{align*}
H^1(\Omega)
&=\{q\in L^2(\Omega):\nabla q\in\bm{L}^2(\Omega)\},\\
\bm{H}(\operatorname{curl},\Omega)
&=\{\bm v\in\bm{L}^2(\Omega):\curl\bm v\in\bm{L}^2(\Omega)\},\\
\bm{H}(\operatorname{div},\Omega)
&=\{\bm v\in\bm{L}^2(\Omega):\diver\bm v\in L^2(\Omega)\}.
\end{align*}
The standard periodic Sobolev spaces are denoted by
\begin{align*}
 \Vper=H^1_{\mathrm{per}}(\Omega), \qquad \Eper=\bm{H}_{\mathrm{per}}(\operatorname{curl},\Omega), \qquad \Fper=\bm{H}_{\mathrm{per}}(\operatorname{div},\Omega).
\end{align*}
We define the zero-mean scalar and vector spaces by
\begin{equation*}
	L^2_0(\Omega)=\Bigl\{q\in L^2(\Omega):
	\int_\Omega q\,\mathrm{d}x=0\Bigr\},\qquad
	\bm{L}^2_0(\Omega)=\Bigl\{\bm v\in\bm{L}^2(\Omega):
	\int_\Omega\bm v\,\mathrm{d}x=\mathbf0\Bigr\}.
\end{equation*}
By periodicity, $\int_\Omega\diver\bm v\,\mathrm{d}x=0$ for every $\bm v\in\Fper$. We therefore set $Q=L^2_0(\Omega)$, which is the pressure space and the endpoint of the periodic complex:
\[
 \Vper\xrightarrow{\gradop}\Eper\xrightarrow{\curlop}\Fper
 \xrightarrow{\divop}Q.
\]
The constant harmonic modes prevent this periodic complex from being exact. To remove these modes, we define
\begin{align*}
 V=\Vper\cap L^2_0(\Omega), \qquad E=\Eper\cap\bm{L}^2_0(\Omega), \qquad F=\Fper\cap\bm{L}^2_0(\Omega).
\end{align*}
The resulting spaces form the exact reduced complex
\begin{equation}
 0\longrightarrow V\xrightarrow{\gradop}E
 \xrightarrow{\curlop}F\xrightarrow{\divop}Q
 \longrightarrow0.
 \label{eq:continuous-exact-complex}
\end{equation}
To pose the MHD problem in these zero-mean spaces, we assume that the prescribed data satisfy
\[
\bm u_0,\bm B_0\in\bm{L}^2_0(\Omega),
\qquad
\bm f(\cdot,t)\in\bm{L}^2_0(\Omega)
\quad\text{for all }t\in[0,T].
\]

\subsection{The dual-field mixed weak formulation}
\label{subsec:weak-formulation}

For the weak formulation and conservation analysis below, we consider sufficiently smooth periodic solutions. Following the dual-field formulation in \cite[Sec.~2.2]{MR4917731}, we use two sets of vector fields, one in $E$ and one in $F$, with corresponding pressure variables in $V$ and $Q$:
\[
\begin{aligned}
 &\bm u_1,\bm B_1,\boldsymbol\omega_1,\bm J_1\in E,
 \qquad p_0\in V,\\
 &\bm u_2,\bm B_2,\boldsymbol\omega_2,\bm J_2\in F,
 \qquad p_3\in Q.
\end{aligned}
\]
The smoothness assumption supplements these space requirements and ensures that all nonlinear terms are well defined; in particular, $\bm u_1\times\bm B_1\in\Eper$, so its curl lies in $\bm{L}^2(\Omega)$.
The continuous dual-field weak formulation is: find the above unknowns such that, for all
\[
\begin{aligned}
 &\bm v_1,\bm z_1,\bm C_1,\bm k_1\in E,
 \qquad q_0\in V,\\
 &\bm v_2,\bm z_2,\bm C_2,\bm k_2\in F,
 \qquad q_3\in Q,
\end{aligned}
\]
the following equations hold:
\begin{subequations}
\label{eq:weak-system-two}
\begin{align}
	\left(\partial_t\bm u_2,\bm v_2\right)
	+\left(\boldsymbol\omega_2\times\bm u_2,\bm v_2\right)
	+Re^{-1}\left(\curl\boldsymbol\omega_1,\bm v_2\right)
	-\left(p_3,\diver\bm v_2\right)
	-s\left(\bm J_2\times\bm B_2,\bm v_2\right)
	&=\left(\bm f,\bm v_2\right),\label{eq:weak-two-momentum}\\
	\left(\boldsymbol\omega_1,\bm z_1\right)
	-\left(\bm u_2,\curl\bm z_1\right)&=0,\label{eq:weak-two-vorticity}\\
	\left(\diver\bm u_2,q_3\right)&=0,\label{eq:weak-two-divergence}\\
	\left(\partial_t\bm B_1,\bm C_1\right)
	+Rm^{-1}\left(\bm J_2,\curl\bm C_1\right)
	-\left(\bm u_2\times\bm B_2,\curl\bm C_1\right)&=0,\label{eq:weak-two-induction}\\
	\left(\bm J_2,\bm k_2\right)
	-\left(\curl\bm B_1,\bm k_2\right)&=0.\label{eq:weak-two-current}
\end{align}
\end{subequations}
and
\begin{subequations}
\label{eq:weak-system-one}
\begin{align}
	\left(\partial_t\bm u_1,\bm v_1\right)
	+\left(\boldsymbol\omega_1\times\bm u_1,\bm v_1\right)
	+Re^{-1}\left(\boldsymbol\omega_2,\curl\bm v_1\right)
	+\left(\nabla p_0,\bm v_1\right)
	-s\left(\bm J_1\times\bm B_1,\bm v_1\right)
	&=\left(\bm f,\bm v_1\right),\label{eq:weak-one-momentum}\\
	\left(\boldsymbol\omega_2,\bm z_2\right)
	-\left(\curl\bm u_1,\bm z_2\right)&=0,\label{eq:weak-one-vorticity}\\
	\left(\bm u_1,\nabla q_0\right)&=0,\label{eq:weak-one-divergence}\\
	\left(\partial_t\bm B_2,\bm C_2\right)
	+Rm^{-1}\left(\curl\bm J_1,\bm C_2\right)
	-\left(\curl(\bm u_1\times\bm B_1),\bm C_2\right)&=0,\label{eq:weak-one-induction}\\
	\left(\bm J_1,\bm k_1\right)
	-\left(\bm B_2,\curl\bm k_1\right)&=0.\label{eq:weak-one-current}
\end{align}
\end{subequations}
Equations~\eqref{eq:weak-one-vorticity} and~\eqref{eq:weak-two-current} imply the following identities in $\bm{L}^2(\Omega)$: $\boldsymbol\omega_2=\curl\bm u_1$ and $\bm J_2=\curl\bm B_1$. Equations~\eqref{eq:weak-two-vorticity} and~\eqref{eq:weak-one-current} impose the corresponding relations for $\boldsymbol\omega_1$ and $\bm J_1$ in weak form, with curl acting on the test functions in $E$ through periodic integration by parts.

\subsection{Divergence constraints and balance laws}
\label{subsec:continuous-conservation}
We derive the divergence constraints and balance laws for the sufficiently
smooth periodic solutions considered in Section~\ref{subsec:weak-formulation}.

\begin{proposition}[Divergence-free constraints]
\label{thm:continuous-divergence}
Assume that the initial magnetic fields satisfy
$\bm B_1(0)\perp\nabla V$ and $\diver\bm B_2(0)=0$.
Then, for all times,
\[
\begin{aligned}
 &\diver\bm u_2=0,\qquad
 \diver\bm B_2=0,\qquad
 \diver\bm J_2=0,\\
 &\bm u_1\perp\nabla V,\qquad
 \bm B_1\perp\nabla V,\qquad
 \bm J_1\perp\nabla V.
\end{aligned}
\]
\end{proposition}
\begin{proof}
The velocity constraints follow directly from
\eqref{eq:weak-two-divergence} and \eqref{eq:weak-one-divergence}.
The curl relations \eqref{eq:weak-two-current} and
\eqref{eq:weak-one-current}, together with the identities
\(\nabla\times\nabla=0\) and \(\nabla\cdot(\nabla\times)=0\),
imply the strong and weak divergence-free properties of the current density.
Applying these identities to the induction equations
\eqref{eq:weak-two-induction} and \eqref{eq:weak-one-induction}
shows that the initial divergence-free conditions on the magnetic field
are preserved in time.
\end{proof}

Define the two energies and the cross helicity by
\begin{equation}
 E_1=\dfrac12\left(\norm{\bm u_1}^2+s\norm{\bm B_2}^2\right),
 \qquad
 E_2=\dfrac12\left(\norm{\bm u_2}^2+s\norm{\bm B_1}^2\right).
 \label{eq:continuous-dual-energies}
\end{equation}
\begin{equation}
	H_c=\dfrac12\left[(\bm u_1,\bm B_1)
	+(\bm u_2,\bm B_2)\right].
	\label{eq:continuous-cross-definition}
\end{equation}

\begin{theorem}[Energy and cross-helicity balances]
\label{thm:continuous-energy-cross}
\label{thm:continuous-energy}
\label{thm:continuous-cross-helicity}
\leavevmode\par\nobreak\smallskip\nobreak
\noindent\textup{\bfseries\makebox[1.8em][l]{(i)}Energy.} The two energies satisfy
\begin{equation}
\begin{aligned}
\dfrac{\mathrm{d}E_1}{\mathrm{d}t}&=-Re^{-1}\norm{\boldsymbol\omega_2}^2
 -sRm^{-1}\norm{\bm J_1}^2+(\bm f,\bm u_1),\\[2pt]
 \dfrac{\mathrm{d}E_2}{\mathrm{d}t}&=-Re^{-1}\norm{\boldsymbol\omega_1}^2
 -sRm^{-1}\norm{\bm J_2}^2+(\bm f,\bm u_2).
\end{aligned}
\label{eq:continuous-energy-laws}
\end{equation}

\par\smallskip
\noindent\textup{\bfseries\makebox[1.8em][l]{(ii)}Cross helicity.} The cross helicity satisfies
\begin{equation}
\begin{aligned}
 \dfrac{\mathrm{d}H_c}{\mathrm{d}t}
 &=-\dfrac12\left(Re^{-1}+Rm^{-1}\right)
 \left[(\boldsymbol\omega_1,\bm J_1)
 +(\boldsymbol\omega_2,\bm J_2)\right]+\dfrac12\left[(\bm f,\bm B_1)
 +(\bm f,\bm B_2)\right].
\end{aligned}
 \label{eq:continuous-cross-law}
\end{equation}
If $\bm f=\mathbf0$ and $Re=Rm=\infty$, then $E_1$, $E_2$,
and $H_c$ are conserved.
\end{theorem}
\begin{proof}
\textup{(i)} \textit{Energy.}
Test \eqref{eq:weak-one-momentum} with $\bm v_1=\bm u_1$ and
\eqref{eq:weak-one-induction} with $\bm C_2=s\bm B_2$.  The divergence constraint and the scalar triple-product identity give
\[
\begin{gathered}
	(\nabla p_0,\bm u_1)=0,\qquad
	(\boldsymbol\omega_1\times\bm u_1,\bm u_1)=0
\end{gathered}
\]
Set $\bm w=\bm u_1\times\bm B_1$, with spatial mean
$\langle\bm w\rangle_\Omega$.
Since $\bm J_1$ has zero mean and constants have zero curl,
\eqref{eq:weak-one-current}, tested with
$\bm k_1=\bm w-\langle\bm w\rangle_\Omega\in E$, gives
$(\curl\bm w,\bm B_2)=(\bm w,\bm J_1)$.
This relation and the scalar triple-product identity show that the Lorentz and induction terms cancel:
\[
\begin{gathered}
 -s(\bm J_1\times\bm B_1,\bm u_1)
 -s(\nabla\times(\bm u_1\times\bm B_1),\bm B_2)=-s(\bm J_1\times\bm B_1,\bm u_1)
 -s(\bm u_1\times\bm B_1,\bm J_1)=0.
\end{gathered}
\]
The relation $\curl\bm u_1=\boldsymbol\omega_2$ from \eqref{eq:weak-one-vorticity} and
\eqref{eq:weak-one-current} with $\bm k_1=\bm J_1$ reduce the diffusion terms
to $Re^{-1}\|\boldsymbol\omega_2\|^2$ and $sRm^{-1}\|\bm J_1\|^2$.
Adding the tested equations proves the first identity in
\eqref{eq:continuous-energy-laws}.

Similarly, test \eqref{eq:weak-two-momentum} with $\bm v_2=\bm u_2$ and
\eqref{eq:weak-two-induction} with $\bm C_1=s\bm B_1$.
Using $\diver\bm u_2=0$ from Proposition~\ref{thm:continuous-divergence},
$\curl\bm B_1=\bm J_2$ from \eqref{eq:weak-two-current}, and the scalar triple-product identity,
we obtain the corresponding zero terms and cancellation:
\[
\begin{gathered}
 (p_3,\diver\bm u_2)=0,\qquad
 (\boldsymbol\omega_2\times\bm u_2,\bm u_2)=0,\qquad
  -s(\bm J_2\times\bm B_2,\bm u_2)
  -s(\bm u_2\times\bm B_2,\nabla\times\bm B_1)=0.
\end{gathered}
\]
Testing \eqref{eq:weak-two-vorticity} with $\bm z_1=\boldsymbol\omega_1$
and using $\curl\bm B_1=\bm J_2$ gives the diffusion terms
$Re^{-1}\|\boldsymbol\omega_1\|^2$ and $sRm^{-1}\|\bm J_2\|^2$.
Adding the tested equations proves the second identity in
\eqref{eq:continuous-energy-laws}.

\smallskip
\noindent\textup{(ii)} \textit{Cross helicity.}
Test \eqref{eq:weak-one-momentum} with $\bm v_1=\bm B_1$ and
\eqref{eq:weak-one-induction} with $\bm C_2=\bm u_2$.
By the same mean-subtraction argument, \eqref{eq:weak-two-vorticity} gives
$(\curl\bm w,\bm u_2)=(\bm w,\boldsymbol\omega_1)$.
The scalar triple-product identity then yields
$(\boldsymbol\omega_1\times\bm u_1,\bm B_1)
-(\curl(\bm u_1\times\bm B_1),\bm u_2)=0$.
The pressure and Lorentz terms vanish as in part~\textup{(i)}.
Using $\curl\bm B_1=\bm J_2$ and \eqref{eq:weak-two-vorticity}
with $\bm z_1=\bm J_1$, we add the tested equations to obtain
\[
 (\partial_t\bm u_1,\bm B_1)+(\partial_t\bm B_2,\bm u_2)
 +Re^{-1}(\boldsymbol\omega_2,\bm J_2)
 +Rm^{-1}(\boldsymbol\omega_1,\bm J_1)=(\bm f,\bm B_1).
\]
Next, test \eqref{eq:weak-two-momentum} with $\bm v_2=\bm B_2$ and
\eqref{eq:weak-two-induction} with $\bm C_1=\bm u_1$.
The pressure, Lorentz, and remaining nonlinear terms vanish or cancel as in part~\textup{(i)}.
Using \eqref{eq:weak-one-current} with $\bm k_1=\boldsymbol\omega_1$ and
$\curl\bm u_1=\boldsymbol\omega_2$, we similarly add the tested equations to obtain
\[
 (\partial_t\bm u_2,\bm B_2)+(\partial_t\bm B_1,\bm u_1)
 +Re^{-1}(\boldsymbol\omega_1,\bm J_1)
 +Rm^{-1}(\boldsymbol\omega_2,\bm J_2)=(\bm f,\bm B_2).
\]
The four time-derivative terms sum to $2\,\mathrm dH_c/\mathrm dt$;
adding these identities proves \eqref{eq:continuous-cross-law}.
The stated conservation laws follow when $Re=Rm=\infty$ and $\bm f=\bm0$.
\end{proof}

To relate these balances to Alfv\'enic alignment, define the averaged energy
and the Els\"asser energies by
\begin{equation}
 \mathcal E=\frac12(E_1+E_2),\qquad
 D_{\pm}=\mathcal E\pm\sqrt{s}\,H_c
 =\frac14\sum_{i=1}^2\|\bm u_i\pm\sqrt{s}\,\bm B_i\|^2.
 \label{eq:continuous-elsasser-energies}
\end{equation}

\begin{corollary}[Els\"asser energies and Alfv\'enic alignment in the ideal case]
\label{thm:continuous-elsasser}
If $Re=Rm=\infty$ and $\bm f=\mathbf0$, then, for either choice of sign,
\begin{equation}
 D_\pm(t)=D_\pm(0).
 \label{eq:continuous-elsasser-conservation}
\end{equation}
In particular, $D_-(0)=0$ implies
$\bm u_i(t)=\sqrt{s}\bm B_i(t)$, whereas $D_+(0)=0$ implies
$\bm u_i(t)=-\sqrt{s}\bm B_i(t)$, for $i=1,2$ and all times.
\end{corollary}
\begin{proof}
Expansion of the squared norms verifies \eqref{eq:continuous-elsasser-energies}.
In the ideal, unforced case, Theorem~\ref{thm:continuous-energy-cross}
gives conservation of $E_1$, $E_2$, and $H_c$. Hence
$D_\pm=\tfrac12(E_1+E_2)\pm\sqrt{s}H_c$ is conserved.
If either initial Els\"asser energy is zero, the corresponding
nonnegative squared norms remain zero, yielding the stated alignment
relations.
\end{proof}

For magnetic helicity, choose sufficiently smooth potentials
$\bm A_1\in E$ and $\bm A_2\in F$ satisfying
\begin{equation}
 \curl\bm A_1=\bm B_2,\qquad
 (\bm A_2,\curl\bm e)=(\bm B_1,\bm e)
 \quad\forall\bm e\in E.
 \label{eq:continuous-potential-relations}
\end{equation}
The periodic Hodge decomposition and the divergence constraints in
Proposition~\ref{thm:continuous-divergence} ensure that these potentials
exist, but they are not unique. Any other zero-mean choices for the same
magnetic fields have the form
$\widetilde{\bm A}_i=\bm A_i+\nabla\phi_i$, with periodic scalar functions
$\phi_i$ of sufficient regularity. By the divergence constraints and
periodic integration by parts, $(\nabla\phi_i,\bm B_i)=0$ for $i=1,2$.
Hence $(\widetilde{\bm A}_i,\bm B_i)=(\bm A_i,\bm B_i)$:
the definitions of $H_{M,1}$ and $H_{M,2}$ below do not depend on the
choice of potentials.
Define the magnetic and fluid helicities by
\[
\begin{aligned}
 H_{M,1}&=(\bm A_1,\bm B_1),&
 H_{M,2}&=(\bm A_2,\bm B_2),\\
 H_{F,1}&=(\bm u_1,\boldsymbol\omega_1),&
 H_{F,2}&=(\bm u_2,\boldsymbol\omega_2).
\end{aligned}
\]

\begin{theorem}[Magnetic- and fluid-helicity balances]
\label{thm:continuous-helicity-balances}
\label{thm:continuous-magnetic-helicity}
\label{thm:continuous-fluid-helicity}
\leavevmode\par\nobreak\smallskip\nobreak
\noindent\textup{\bfseries\makebox[1.8em][l]{(i)}Magnetic helicity.} The two magnetic helicities coincide.
Their common value $H_M=H_{M,1}=H_{M,2}$ satisfies
\begin{equation}
		\frac{\mathrm{d}H_M}{\mathrm{d}t}=
		-Rm^{-1}\left[
		(\bm J_1,\bm B_1)+(\bm J_2,\bm B_2)
		\right].
		\label{eq:continuous-magnetic-balance}
	\end{equation}
In particular, $H_M$ is conserved if $Rm=\infty$.

\par\smallskip
\noindent\textup{\bfseries\makebox[1.8em][l]{(ii)}Fluid helicity.} The two fluid helicities coincide.
Their common value $H_F=H_{F,1}=H_{F,2}$ satisfies
\begin{equation}
\begin{aligned}
 \frac{\mathrm dH_F}{\mathrm dt}
 &=-2Re^{-1}(\boldsymbol\omega_2,\curl\boldsymbol\omega_1)+s\bigl[
 (\bm J_1\times\bm B_1,\boldsymbol\omega_1)
 +(\bm J_2\times\bm B_2,\boldsymbol\omega_2)\bigr]+(\bm f,\boldsymbol\omega_1)+(\bm f,\boldsymbol\omega_2).
\end{aligned}
\label{eq:continuous-fluid-balance}
\end{equation}
In particular, $H_F$ is conserved in the unforced inviscid
hydrodynamic case $\bm B_1=\bm B_2=\mathbf0$,
$\bm f=\mathbf0$, and $Re=\infty$.
\end{theorem}
\begin{proof}
\textup{(i)} \textit{Magnetic helicity.}
Testing \eqref{eq:continuous-potential-relations} with
$\bm e=\bm A_1$ and $\bm e=\partial_t\bm A_1$ gives
$H_{M,1}=H_{M,2}$ and
\[
 \frac{\mathrm dH_M}{\mathrm dt}
 =(\partial_t\bm B_2,\bm A_2)+(\partial_t\bm B_1,\bm A_1).
\]
Test \eqref{eq:weak-one-induction} with $\bm C_2=\bm A_2$ and
\eqref{eq:weak-two-induction} with $\bm C_1=\bm A_1$.
For $\bm w=\bm u_1\times\bm B_1$, apply the weak potential relation to
$\bm e=\bm w-\langle\bm w\rangle_\Omega$, using the zero mean of $\bm B_1$.
Together with $\curl\bm A_1=\bm B_2$, this gives
\[
 (\curl\bm w,\bm A_2)=(\bm w,\bm B_1)=0,
 \qquad (\bm u_2\times\bm B_2,\nabla\times\bm A_1)=0.
\]
The potential relation with $\bm e=\bm J_1$ gives
$(\curl\bm J_1,\bm A_2)=(\bm J_1,\bm B_1)$.
Thus the resistive terms are $Rm^{-1}(\bm J_1,\bm B_1)$ and
$Rm^{-1}(\bm J_2,\bm B_2)$, and adding the tested equations proves
\eqref{eq:continuous-magnetic-balance}.
When $Rm=\infty$, its right-hand side vanishes.

\smallskip
\noindent\textup{(ii)} \textit{Fluid helicity.}
Testing \eqref{eq:weak-two-vorticity} with $\bm z_1=\bm u_1$ and
$\bm z_1=\partial_t\bm u_1$, and using
$\curl\bm u_1=\boldsymbol\omega_2$, gives $H_{F,1}=H_{F,2}$ and
\[
 \frac{\mathrm dH_F}{\mathrm dt}
 =(\partial_t\bm u_1,\boldsymbol\omega_1)
 +(\partial_t\bm u_2,\boldsymbol\omega_2).
\]
Test \eqref{eq:weak-one-momentum} with $\bm v_1=\boldsymbol\omega_1$ and
\eqref{eq:weak-two-momentum} with $\bm v_2=\boldsymbol\omega_2$.
The weak vorticity relation \eqref{eq:weak-two-vorticity} with $\bm z_1=\nabla p_0$,
$\diver\boldsymbol\omega_2=\diver\curl\bm u_1=0$, and the triple-product identity give
\[
\begin{gathered}
 (\nabla p_0,\boldsymbol\omega_1)=0,\qquad
 (p_3,\diver\boldsymbol\omega_2)=0,\qquad (\boldsymbol\omega_1\times\bm u_1,\boldsymbol\omega_1)=0,
 \qquad
 (\boldsymbol\omega_2\times\bm u_2,\boldsymbol\omega_2)=0.
\end{gathered}
\]
The two viscous terms are both
$Re^{-1}(\boldsymbol\omega_2,\curl\boldsymbol\omega_1)$.
Adding the tested equations therefore proves \eqref{eq:continuous-fluid-balance}.
Its right-hand side vanishes in the stated unforced inviscid hydrodynamic case.
\end{proof}

\section{A staggered dual-field scheme}
\label{sec:discrete-scheme}

\subsection{Periodic finite element de Rham complex}
\label{subsec:periodic-fe-complex}

On a shape-regular periodic tetrahedral mesh, let $k\ge1$. We use the zero-mean periodic finite element spaces $V_h$, $E_h$, $F_h$, and $Q_h$: continuous Lagrange polynomials of degree $k$, first-family N\'ed\'elec elements of order $k$, Raviart--Thomas elements of degree $k-1$, and discontinuous polynomials of degree $k-1$, respectively. Periodic degrees of freedom are identified and the scalar and vector means are removed, so $V_h\subset V$, $E_h\subset E$, $F_h\subset F$, and $Q_h\subset Q$. These spaces form the exact reduced complex
\begin{equation}
 0\longrightarrow V_h\xrightarrow{\gradop}E_h
 \xrightarrow{\curlop}F_h\xrightarrow{\divop}Q_h
 \longrightarrow0.
 \label{eq:discrete-exact-complex}
\end{equation}
The finite element de Rham complexes and commuting projections are described in \cite{MR2269741,MR2009375}; arbitrary-order Cartesian and simplicial constructions are given in \cite{MR3670721}. The weak curl $\curlh:\Fh\to\Eh$ is the $L^2$ adjoint of the strong discrete
curl from $\Eh$ to $\Fh$:
\begin{equation}
	\left(\curlh \bm{c}_h,\bm{e}_h\right)
	=\left(\bm{c}_h,\curl \bm{e}_h\right),
	\qquad \forall \bm{e}_h\in\Eh.
	\label{eq:weakcurl}
\end{equation}

Let $\Pi_h^E:[L^2(\Omega)]^3\to\Eh$ denote the $L^2$-orthogonal projection, defined by
\begin{equation*}
	(\Pi_h^E\bm g,\bm e_h)=(\bm g,\bm e_h),
	\qquad\forall\bm e_h\in\Eh.
\end{equation*}
Since $\curlh\bm C_h\in\Eh$ for every $\bm C_h\in\Fh$, it follows that
\begin{equation}
	(\Pi_h^E\bm g,\curlh\bm C_h)
	=(\bm g,\curlh\bm C_h),
	\qquad
	\forall\bm g\in[L^2(\Omega)]^3.
	\label{eq:projection-weak-curl}
\end{equation}
Thus, when a nonlinear product does not belong to $\Eh$, its $L^2$ projection onto $\Eh$ may be inserted or omitted in a pairing with $\curlh\bm C_h$ without changing the resulting discrete functional. The spatially discrete formulation uses these finite element spaces in place of the corresponding continuous Sobolev spaces.

\subsection{Staggered integer-step and half-integer-step updates}
\label{subsec:staggered-updates}

Let $N$ be the number of full time steps and set $\dt=T/N$.
The integer grid is $t_n=n\dt$ for $0\le n\le N$ and ends at
$t_N=T$. The half-integer grid is
$t_{n+1/2}=(n+1/2)\dt$ for $0\le n\le N$ and includes the
additional point $t_{N+1/2}=T+\dt/2$. The integer-time fields are
\[
 (\bm u_2^n,\boldsymbol\omega_1^n,\bm B_1^n,\bm J_2^n),
\]
and the half-integer-time fields are
\[
 (\bm u_1^{n+1/2},\boldsymbol\omega_2^{n+1/2},
  \bm B_2^{n+1/2},\bm J_1^{n+1/2}).
\]
The pressures are determined at the midpoints of the update
intervals: $p_3^{n+1/2}$ $(0\le n\le N-1)$ in the integer-step update and $p_0^n$ $(1\le n\le N)$ in the half-integer-step update.
The magnetic potentials $\bm A_1^{n+1/2}$ and
$\bm A_2^n$ are reconstructed from
$\bm B_2^{n+1/2}$ and $\bm B_1^n$, respectively.
They are used in the magnetic-helicity analysis and are
not included among the unknowns in the time-stepping systems.
For a half-integer-step variable $b$ and an integer-step variable $a$, define
\begin{alignat*}{2}
	\bar b^n
	&=\tfrac12\bigl(b^{n-1/2}+b^{n+1/2}\bigr),
	&\qquad \delta_t b^n
	&=\bigl(b^{n+1/2}-b^{n-1/2}\bigr)/\dt,\\
	\bar a^{n+1/2}
	&=\tfrac12\bigl(a^n+a^{n+1}\bigr),
	&\qquad \delta_t a^{n+1/2}
	&=\bigl(a^{n+1}-a^n\bigr)/\dt.
\end{alignat*}

At $t=0$, in addition to the integer-step fields
$\bm u_2^0$ and $\bm B_1^0$, introduce the auxiliary values $\bm u_1^0$ and $\bm B_2^0$. They satisfy
\[
\bm B_1^0\perp\nabla\Vh, \qquad \bm u_1^0\perp\nabla\Vh,
\qquad
\diver\bm u_2^0=\diver\bm B_2^0=0.
\]
The corresponding initial vorticities and currents are defined by
\[
\boldsymbol\omega_1^0=\curlh\bm u_2^0,
\qquad
\bm J_2^0=\curl\bm B_1^0,
\qquad
\boldsymbol\omega_2^0=\curl\bm u_1^0,
\qquad
\bm J_1^0=\curlh\bm B_2^0.
\]

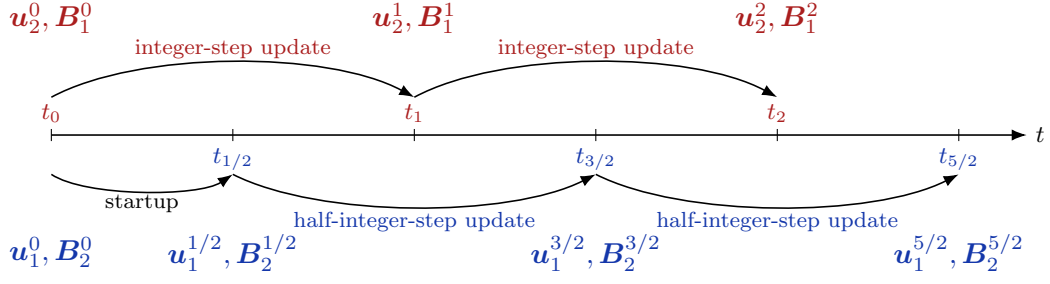
\begin{figure}[htbp]
\centering
\begingroup
\definecolor{integerstep}{RGB}{170,30,30}
\definecolor{halfstep}{RGB}{20,55,170}
\begin{tikzpicture}[x=1cm,y=1cm,>=Latex,font=\normalsize,
  advance/.style={-{Latex[length=2mm]},line width=0.6pt}]
  \draw[advance,black] (1,0) -- (13.9,0) node[right] {$t$};
  \foreach \x/\time in {1/{0},5.8/{1},10.6/{2}} {
    \draw (\x,0.08) -- (\x,-0.08);
    \node[above,inner sep=1pt,text=integerstep] at (\x,0.13) {$t_{\time}$};
  }
  \foreach \x/\time in {3.4/{1/2},8.2/{3/2},13/{5/2}} {
    \draw (\x,0.08) -- (\x,-0.08);
    \node[below,inner sep=1pt,text=halfstep] at (\x,-0.13) {$t_{\time}$};
  }
  \node[align=center,font=\large,text=integerstep] at (1,1.55) {$\bm u_2^0,\bm B_1^0$};
  \node[align=center,font=\large,text=halfstep] at (1,-1.55)
    {$\bm u_1^0,\bm B_2^0$};
  \node[align=center,font=\large,text=halfstep] at (3.4,-1.55) {$\bm u_1^{1/2},\bm B_2^{1/2}$};
  \node[align=center,font=\large,text=integerstep] at (5.8,1.55) {$\bm u_2^1,\bm B_1^1$};
  \node[align=center,font=\large,text=halfstep] at (8.2,-1.55) {$\bm u_1^{3/2},\bm B_2^{3/2}$};
  \node[align=center,font=\large,text=integerstep] at (10.6,1.55) {$\bm u_2^2,\bm B_1^2$};
  \node[align=center,font=\large,text=halfstep] at (13,-1.55) {$\bm u_1^{5/2},\bm B_2^{5/2}$};
  \draw[advance,black] (1,-0.52) .. controls (1.45,-0.82) and (2.95,-0.82) .. (3.4,-0.52);
  \node[font=\small] at (2.2,-0.90) {startup};
  \draw[advance,black] (1,0.50) .. controls (2.15,1.12) and (4.65,1.12) .. (5.8,0.50);
  \node[font=\small,text=integerstep] at (3.4,1.14) {integer-step update};
  \draw[advance,black] (3.4,-0.52) .. controls (4.55,-1.10) and (7.05,-1.10) .. (8.2,-0.52);
  \node[font=\small,text=halfstep] at (5.8,-1.14) {half-integer-step update};
  \draw[advance,black] (5.8,0.50) .. controls (6.95,1.12) and (9.45,1.12) .. (10.6,0.50);
  \node[font=\small,text=integerstep] at (8.2,1.14) {integer-step update};
  \draw[advance,black] (8.2,-0.52) .. controls (9.35,-1.10) and (11.85,-1.10) .. (13,-0.52);
  \node[font=\small,text=halfstep] at (10.6,-1.14) {half-integer-step update};
\end{tikzpicture}
\endgroup
\caption{Staggered storage and updates for the primary velocity and magnetic fields. Integer times and states are shown in red above the axis, and half-integer times and states in blue below it. Black arrows connect successive states on each grid; the first lower arrow represents the half-length startup solve. Thereafter the integer-step and half-integer-step solves alternate, each using the most recently available state on the other grid.}
\label{fig:staggered-field-timeline}
\end{figure}

The first half-integer-step fields at $t=\dt/2$ are computed using Algorithm~\ref{alg:half-update}. In this startup solve, the time increment is $\dt/2$, the evolving fields are averaged between $t=0$ and $t=\dt/2$, and $\boldsymbol\omega_1^0$ and $\bm B_1^0$ are held fixed. If the body force depends on time, it is evaluated at $t=\dt/4$. The resulting fields satisfy
\[
\boldsymbol\omega_2^{1/2}
=\curl\bm u_1^{1/2},
\qquad
\bm J_1^{1/2}
=\curlh\bm B_2^{1/2},
\qquad
\bm u_1^{1/2}\perp\nabla\Vh,
\qquad
\diver\bm B_2^{1/2}=0.
\]
After this startup solve, each time step applies Algorithm~\ref{alg:integer-update} followed by Algorithm~\ref{alg:half-update}, as shown in Figure~\ref{fig:staggered-field-timeline}. The same discrete curl relations and divergence-free constraints hold at all subsequent time levels.

\begin{algorithm}[Integer-step update]
\label{alg:integer-update}
For $0\le n\le N-1$, find
$(\bm{u}_2^{n+1},\boldsymbol{\omega}_1^{n+1},p_3^{n+1/2},\bm{B}_1^{n+1},\bm{J}_2^{n+1})$ in $\Fh\times\Eh\times\Qh\times\Eh\times\Fh$ such that
\begingroup

\begin{subequations}
\label{eq:integer-update-system}
\begin{align}
	\bigl(\delta_t \bm{u}_2^{n+1/2},\bm{v}_2\bigr)
	+\bigl(\boldsymbol{\omega}_2^{n+1/2}\times\overline{\bm{u}}_2^{n+1/2},\bm{v}_2\bigr)+Re^{-1}\bigl(\curl\overline{\boldsymbol{\omega}}_1^{n+1/2},\bm{v}_2\bigr)
	\nonumber\\
	-\bigl(p_3^{n+1/2},\diver \bm{v}_2\bigr)
	-s\bigl(\overline{\bm{J}}_2^{n+1/2}\times \bm{B}_2^{n+1/2},\bm{v}_2\bigr)
	&=\bigl(\bm{f}^{n+1/2},\bm{v}_2\bigr),\label{eq:int-u}\\
	\bigl(\boldsymbol{\omega}_1^{n+1},\bm{z}_1\bigr)-\bigl(\bm{u}_2^{n+1},\curl \bm{z}_1\bigr)
	&=0,\label{eq:int-omega}\\
	\bigl(\diver \bm{u}_2^{n+1},q_3\bigr)&=0,\label{eq:int-p}\\
	\bigl(\delta_t \bm{B}_1^{n+1/2},\bm{C}_1\bigr)
	+Rm^{-1}\bigl(\overline{\bm{J}}_2^{n+1/2},\curl \bm{C}_1\bigr)
	-\bigl(\overline{\bm{u}}_2^{n+1/2}\times \bm{B}_2^{n+1/2},\curl \bm{C}_1\bigr)&=0,
	\label{eq:int-B}\\
	\bigl(\bm{J}_2^{n+1},\bm{k}_2\bigr)-\bigl(\curl \bm{B}_1^{n+1},\bm{k}_2\bigr)&=0,
	\label{eq:int-J}
\end{align}
\end{subequations}
\endgroup
for any $(\bm{v}_2,\bm{z}_1,q_3,\bm{C}_1,\bm{k}_2)$ in $\Fh\times\Eh\times\Qh\times\Eh\times\Fh$.
\end{algorithm}

\begin{algorithm}[Half-integer-step update]
\label{alg:half-update}

For $1\le n\le N$, find
$(\bm{u}_1^{n+1/2},\boldsymbol{\omega}_2^{n+1/2},p_0^n,\bm{B}_2^{n+1/2},\bm{J}_1^{n+1/2})$ in $\Eh\times\Fh\times\Vh\times\Fh\times\Eh$ such that
\begingroup

\begin{subequations}
\label{eq:half-update-system}
\begin{align}
	\bigl(\delta_t \bm{u}_1^n,\bm{v}_1\bigr)+\bigl(\boldsymbol{\omega}_1^n\times\overline{\bm{u}}_1^n,\bm{v}_1\bigr)+Re^{-1}\bigl(\overline{\boldsymbol{\omega}}_2^n,\curl \bm{v}_1\bigr)
	+\bigl(\nabla p_0^n,\bm{v}_1\bigr)
	-s\bigl(\overline{\bm{J}}_1^n\times \bm{B}_1^n,\bm{v}_1\bigr)
	&=\bigl(\bm{f}^n,\bm{v}_1\bigr),\label{eq:half-u}\\
	\bigl(\boldsymbol{\omega}_2^{n+1/2},\bm{z}_2\bigr)-\bigl(\curl \bm{u}_1^{n+1/2},\bm{z}_2\bigr)
	&=0,\label{eq:half-omega}\\
	\bigl(\bm{u}_1^{n+1/2},\nabla q_0\bigr)&=0,\label{eq:half-p}\\
	\bigl(\delta_t \bm{B}_2^n,\bm{C}_2\bigr)+Rm^{-1}\bigl(\overline{\bm{J}}_1^n,\curlh \bm{C}_2\bigr)
	-\bigl(\overline{\bm{u}}_1^n\times \bm{B}_1^n,\curlh \bm{C}_2\bigr)&=0,
	\label{eq:half-B}\\
	\bigl(\bm{J}_1^{n+1/2},\bm{k}_1\bigr)-\bigl(\bm{B}_2^{n+1/2},\curl \bm{k}_1\bigr)
	&=0.\label{eq:half-J}
\end{align}
\end{subequations}
\endgroup
for any $(\bm{v}_1,\bm{z}_2,q_0,\bm{C}_2,\bm{k}_1)$ in $\Eh\times\Fh\times\Vh\times\Fh\times\Eh$.
\end{algorithm}

\begin{remark}
	By the adjoint relation \eqref{eq:weakcurl} and the projection
	property \eqref{eq:projection-weak-curl},
	\eqref{eq:half-B} is equivalent to
	\begin{equation}
		\left(
		\delta_t\bm B_2^n
		+\curl\left[
		Rm^{-1}\overline{\bm J}_1^n
		-\Pi_h^E\bigl(\overline{\bm u}_1^n\times\bm B_1^n\bigr)
		\right],\bm C_2
		\right)=0,
		\qquad \forall\bm C_2\in F_h.
		\label{eq:half-B-projected-curl}
	\end{equation}
\end{remark}

\begin{remark}[Numerical implementation]
	The nonlinear induction term in \eqref{eq:half-B} is
	implemented either using the assembled weak-curl matrix
	or through a periodic edge-element mass-matrix projection followed by the discrete curl.
\end{remark}

Table~\ref{tab:staggered-pairing} compares the velocity--magnetic unknown pairs of the proposed method and the Mao--Xi scheme \cite{MR4917731}. The proposed choice of unknown pairs, together with the associated induction and Lorentz couplings, provides the nonlinear cancellations needed for an additional cross-helicity balance while retaining the energy balances. Each substep remains linear because the fields used as coupling coefficients have already been computed in the other staggered update. The cancellations for energy and cross helicity are established in Theorem~\ref{thm:discrete-energy-cross}.

\subsection{Discrete divergence constraints and balance laws}
\label{subsec:discrete-properties}
We organize the fully discrete results in the same order as the continuous
identities in Section~\ref{subsec:continuous-conservation}: divergence
constraints; energy and cross-helicity balances and their relation to Alfv\'enic alignment;
and magnetic and fluid helicity. The magnetic-helicity law also yields
positive energy lower bounds. Table~\ref{tab:results-guide} provides a guide
to the continuous and discrete conservation results, their sufficient
conditions, and the corresponding numerical tests.

\begin{table}[!htbp]
\centering

\caption{Continuous and discrete conservation laws and their numerical tests.}
\label{tab:results-guide}
\small
\setlength{\tabcolsep}{3pt}
\renewcommand{\arraystretch}{1.15}
\begin{tabularx}{\linewidth}{@{}>{\raggedright\arraybackslash}p{0.17\linewidth}
>{\raggedright\arraybackslash}X
>{\raggedright\arraybackslash}p{0.18\linewidth}
>{\raggedright\arraybackslash}p{0.18\linewidth}
>{\raggedright\arraybackslash}p{0.18\linewidth}@{}}
\toprule
Conserved quantity & Conditions & Continuous result & Discrete result & Numerical test\\
\midrule
Energy
& $Re=Rm=\infty$, $\bm f=\mathbf0$
& Theorem~\ref{thm:continuous-energy-cross}(i)
& Theorem~\ref{thm:discrete-energy-cross}(i)
& Section~\ref{subsubsec:energy-cross-divergence-test}\\
\addlinespace
Cross helicity
& $Re=Rm=\infty$, $\bm f=\mathbf0$
& Theorem~\ref{thm:continuous-energy-cross}(ii)
& Theorem~\ref{thm:discrete-energy-cross}(ii)
& Sections~\ref{subsubsec:aligned-ideal-test} and \ref{subsubsec:energy-cross-divergence-test}\\
\addlinespace
Magnetic helicity
& $Rm=\infty$
& Theorem~\ref{thm:continuous-helicity-balances}(i)
& Theorem~\ref{thm:discrete-helicity-balances}(i)
& Section~\ref{subsubsec:magnetic-helicity-test}\\
\addlinespace
Fluid helicity
& $Re=\infty$, $\bm f=\mathbf0$, $\bm B_1=\bm B_2=\mathbf0$
& Theorem~\ref{thm:continuous-helicity-balances}(ii)
& Theorem~\ref{thm:discrete-helicity-balances}(ii)
& Section~\ref{subsubsec:fluid-helicity-test}\\
\bottomrule
\end{tabularx}
\end{table}

\subsubsection{Divergence-free constraints}
\label{subsubsec:divergence-preservation}

\begin{proposition}[Divergence-free constraints]
\label{thm:discrete-divergence-preservation}
With the initialization and startup described in Section~\ref{subsec:staggered-updates}, the following constraints hold for $0\le n\le N$:
\[
\begin{aligned}
 &\diver\bm u_2^n=0,\qquad
 \diver\bm B_2^{n+1/2}=0,\qquad
 \diver\bm J_2^n=0,\\
 &\bm u_1^{n+1/2}\perp\nabla V_h,\qquad
 \bm B_1^n\perp\nabla V_h,\qquad
 \bm J_1^{n+1/2}\perp\nabla V_h.
\end{aligned}
\]
\end{proposition}
\begin{proof}
At each updated time level, the velocity constraints follow directly from \eqref{eq:int-p} and \eqref{eq:half-p}.
	Gradient tests in \eqref{eq:int-B} and \eqref{eq:half-J}, together with $\nabla\times\nabla=0$, propagate the weak
	magnetic constraint and give the weak divergence-free property of the current density. The strong divergence-free properties of the magnetic field and current density follow from \eqref{eq:half-B-projected-curl}
	and \eqref{eq:int-J}, respectively, using
	$\nabla\cdot(\nabla\times)=0$. The initialized integer-time fields and the startup half-integer-step fields supply the constraints at $n=0$.
\end{proof}

\subsubsection{Energy and cross-helicity balances}
\label{subsubsec:energy-conservation}
\label{subsubsec:cross-helicity}

For $0\le n\le N$, define $E_1^{n+1/2}$ and $E_2^n$ by
\begin{align}
	E_1^{n+1/2}&=\dfrac{1}{2}\bigl\|\bm{u}_1^{n+1/2}\bigr\|^2
	+\dfrac{s}{2}\bigl\|\bm{B}_2^{n+1/2}\bigr\|^2,\label{eq:E1}\\
	E_2^n&=\dfrac{1}{2}\norm{\bm{u}_2^n}^2+\dfrac{s}{2}\norm{\bm{B}_1^n}^2.
	\label{eq:E2}
\end{align}
For $1\le n\le N$, define the integer-step cross helicity
\begin{equation}
	H_c^n=\dfrac12\left(\left(\overline{\bm{u}}_1^n,\bm{B}_1^n\right)
	+\left(\bm{u}_2^n,\overline{\bm{B}}_2^n\right)\right).
	\label{eq:Hc}
\end{equation}
No half-integer-step cross helicity is introduced.

\begin{theorem}[Discrete energy and cross-helicity balances]
\label{thm:discrete-energy-cross}
\label{thm:discrete-energy-laws}
\label{thm:cross-helicity}
\leavevmode\par\nobreak\smallskip\nobreak
\noindent\textup{\bfseries\makebox[1.8em][l]{(i)}Energy.} For $1\le n\le N$,
\begin{align}
 E_1^{n+1/2}-E_1^{n-1/2}
 &=-\dt Re^{-1}\bigl\|\overline{\boldsymbol{\omega}}_2^n\bigr\|^2
 -s\dt Rm^{-1}\bigl\|\overline{\bm{J}}_1^n\bigr\|^2
 +\dt\bigl(\bm{f}^n,\overline{\bm{u}}_1^n\bigr),
 \label{eq:energy1-law}\\[2pt]
 E_2^n-E_2^{n-1}
 &=-\dt Re^{-1}\bigl\|\overline{\boldsymbol{\omega}}_1^{n-1/2}\bigr\|^2
 -s\dt Rm^{-1}\bigl\|\overline{\bm{J}}_2^{n-1/2}\bigr\|^2
 +\dt\bigl(\bm{f}^{n-1/2},\overline{\bm{u}}_2^{n-1/2}\bigr).
 \label{eq:energy2-law}
\end{align}

\par\smallskip
\noindent\textup{\bfseries\makebox[1.8em][l]{(ii)}Cross helicity.} For $1\le n\le N-1$,
\begin{equation}
\begin{aligned}
 H_c^{n+1}-H_c^n
 &=-\frac{\dt}{4}\bigl(
 \mathcal D_{H,1}^n+2\mathcal D_{H,2}^{n+1/2}
 +\mathcal D_{H,1}^{n+1}\bigr)+\frac{\dt}{4}\bigl(
 \mathcal F_{H,1}^n+2\mathcal F_{H,2}^{n+1/2}
 +\mathcal F_{H,1}^{n+1}\bigr),
\end{aligned}
\label{eq:cross-helicity-law}
\end{equation}
where, for $1\le n\le N$, $\mathcal D_{H,i}$ and $\mathcal F_{H,i}$ denote the diffusion and forcing contributions to cross helicity, respectively:
\begin{align}
		\mathcal D_{H,1}^n
		&=Re^{-1}\bigl(\overline{\boldsymbol\omega}_2^n,\bm J_2^n\bigr)
		+Rm^{-1}\bigl(\boldsymbol\omega_1^n,\overline{\bm J}_1^n\bigr),
		\label{eq:cross-dissipation-one}\\
		\mathcal D_{H,2}^{n-1/2}
		&=Re^{-1}\bigl(\overline{\boldsymbol\omega}_1^{n-1/2},
		\bm J_1^{n-1/2}\bigr)
		+Rm^{-1}\bigl(\boldsymbol\omega_2^{n-1/2},
		\overline{\bm J}_2^{n-1/2}\bigr),
		\label{eq:cross-dissipation-two}\\
		\mathcal F_{H,1}^n&=\bigl(\bm f^n,\bm B_1^n\bigr), \qquad
		\mathcal F_{H,2}^{n-1/2}
		=\bigl(\bm f^{n-1/2},\bm B_2^{n-1/2}\bigr).
		\label{eq:cross-forcing}
	\end{align}
The subscript $H$ refers to cross helicity $H_c$; indices $1$ and $2$
label the half-integer-step and integer-step updates, respectively.
If $\bm f=\mathbf0$ and $Re=Rm=\infty$, both energy
sequences and $H_c^n$ are conserved on their respective index ranges.
\end{theorem}
\begin{proof}
\textup{(i)} \textit{Energy.}
Test \eqref{eq:half-u} with $\bm v_1=\overline{\bm u}_1^n$ and
\eqref{eq:half-B} with $\bm C_2=s\overline{\bm B}_2^n$.
The time terms give
\[
 \dt\bigl[(\delta_t\bm u_1^n,\overline{\bm u}_1^n)
 +s(\delta_t\bm B_2^n,\overline{\bm B}_2^n)\bigr]
 =E_1^{n+1/2}-E_1^{n-1/2}.
\]
Averaging \eqref{eq:half-J} and \eqref{eq:half-omega} gives 
$\curlh\overline{\bm B}_2^n=\overline{\bm J}_1^n$ and
$\curl\overline{\bm u}_1^n=\overline{\boldsymbol\omega}_2^n$, respectively.
Using these curl relations, the divergence constraint, and the scalar triple-product identity, we obtain
\[
\begin{gathered}
 (\nabla p_0^n,\overline{\bm u}_1^n)=0,
 \qquad
 (\boldsymbol\omega_1^n\times\overline{\bm u}_1^n,\overline{\bm u}_1^n)=0,\\
 -s(\overline{\bm J}_1^n\times\bm B_1^n,\overline{\bm u}_1^n)
 -s(\overline{\bm u}_1^n\times\bm B_1^n,\curlh\overline{\bm B}_2^n)=0.
\end{gathered}
\]
The curl relations reduce the two diffusion terms to
$Re^{-1}\|\overline{\boldsymbol\omega}_2^n\|^2$ and
$sRm^{-1}\|\overline{\bm J}_1^n\|^2$, respectively.
Adding the tested equations and multiplying by $\dt$ proves \eqref{eq:energy1-law}.

For the integer-step identity, test \eqref{eq:int-u} and \eqref{eq:int-B} at step $n-1$
with $\bm v_2=\overline{\bm u}_2^{n-1/2}$ and
$\bm C_1=s\overline{\bm B}_1^{n-1/2}$, respectively.
The time terms similarly give
\[
 \dt\bigl[(\delta_t\bm u_2^{n-1/2},\overline{\bm u}_2^{n-1/2})
 +s(\delta_t\bm B_1^{n-1/2},\overline{\bm B}_1^{n-1/2})\bigr]
 =E_2^n-E_2^{n-1}.
\]
Averaging \eqref{eq:int-J} and \eqref{eq:int-omega} gives
$\curl\overline{\bm B}_1^{n-1/2}=\overline{\bm J}_2^{n-1/2}$ and
$\curlh\overline{\bm u}_2^{n-1/2}=\overline{\boldsymbol\omega}_1^{n-1/2}$, respectively.
Using these curl relations, the divergence constraint, and the scalar triple-product identity, we similarly obtain
\[
\begin{gathered}
 (p_3^{n-1/2},\diver\overline{\bm u}_2^{n-1/2})=0,
 \qquad
 (\boldsymbol\omega_2^{n-1/2}\times\overline{\bm u}_2^{n-1/2},\overline{\bm u}_2^{n-1/2})=0,\\
 -s(\overline{\bm J}_2^{n-1/2}\times\bm B_2^{n-1/2},\overline{\bm u}_2^{n-1/2})
 -s(\overline{\bm u}_2^{n-1/2}\times\bm B_2^{n-1/2},\curl\overline{\bm B}_1^{n-1/2})=0.
\end{gathered}
\]
The curl relations reduce the two diffusion terms to
$Re^{-1}\|\overline{\boldsymbol\omega}_1^{n-1/2}\|^2$ and
$sRm^{-1}\|\overline{\bm J}_2^{n-1/2}\|^2$, respectively.
Adding the tested equations and multiplying by $\dt$ proves \eqref{eq:energy2-law}.
If $Re=Rm=\infty$ and $\bm f=\bm0$, the dissipative and forcing terms vanish, so
\begin{equation}
 E_1^{n+1/2}=E_1^{n-1/2},\qquad E_2^n=E_2^{n-1}.
 \label{eq:energy-conservation-explicit}
\end{equation}

\smallskip
\noindent\textup{(ii)} \textit{Cross helicity.}
Test \eqref{eq:half-u} with $\bm v_1=\bm B_1^n$ and
\eqref{eq:half-B} with $\bm C_2=\bm u_2^n$.
Equations~\eqref{eq:int-J} and \eqref{eq:int-omega} give
$\curl\bm B_1^n=\bm J_2^n$ and
$\curlh\bm u_2^n=\boldsymbol\omega_1^n$, respectively.
The divergence constraint and the scalar triple-product identities give
\[
\begin{gathered}
 (\nabla p_0^n,\bm B_1^n)=0,
 \qquad (\overline{\bm J}_1^n\times\bm B_1^n,\bm B_1^n)=0,\\
 (\boldsymbol\omega_1^n\times\overline{\bm u}_1^n,\bm B_1^n)
 -(\overline{\bm u}_1^n\times\bm B_1^n,\curlh\bm u_2^n)=0.
\end{gathered}
\]
Adding the tested equations and multiplying by $\dt$ therefore gives
\begin{equation}
 (\bm u_1^{n+1/2}-\bm u_1^{n-1/2},\bm B_1^n)
 +(\bm u_2^n,\bm B_2^{n+1/2}-\bm B_2^{n-1/2})
 +\dt\mathcal D_{H,1}^n
 =\dt\mathcal F_{H,1}^n.
 \label{eq:Hc-half-bridge}
\end{equation}

Next, test \eqref{eq:int-u} with $\bm v_2=\bm B_2^{n+1/2}$ and
\eqref{eq:int-B} with $\bm C_1=\bm u_1^{n+1/2}$.
Equations~\eqref{eq:half-J} and \eqref{eq:half-omega} give
$\curlh\bm B_2^{n+1/2}=\bm J_1^{n+1/2}$ and
$\curl\bm u_1^{n+1/2}=\boldsymbol\omega_2^{n+1/2}$, respectively.
The pressure and Lorentz terms vanish, and the remaining nonlinear terms
cancel as above. Using the adjoint relation, adding the tested equations,
and multiplying by $\dt$ gives
\begin{equation}
 (\bm u_1^{n+1/2},\bm B_1^{n+1}-\bm B_1^n)
 +(\bm u_2^{n+1}-\bm u_2^n,\bm B_2^{n+1/2})
 +\dt\mathcal D_{H,2}^{n+1/2}
 =\dt\mathcal F_{H,2}^{n+1/2}.
 \label{eq:Hc-int-bridge}
\end{equation}
Add \eqref{eq:Hc-half-bridge} at steps $n$ and $n+1$ to twice
\eqref{eq:Hc-int-bridge}. By \eqref{eq:Hc}, the inner-product increments
combine to $4(H_c^{n+1}-H_c^n)$; dividing the resulting identity by four
proves \eqref{eq:cross-helicity-law}.
In the ideal unforced case, all $\mathcal D_{H,i}$ and $\mathcal F_{H,i}$ vanish,
so $H_c^{n+1}=H_c^n$.
\end{proof}

\begin{remark}[Cumulative cross-helicity balance]
\label{rem:cross-helicity-cumulative}
Summing \eqref{eq:cross-helicity-law} gives the cumulative form used for
the balance diagnostic in the numerical experiments: for $2\le\ell\le N$,
\begin{equation}
\begin{aligned}
 H_c^\ell-H_c^1
 &=-\frac{\dt}{4}\sum_{n=1}^{\ell-1}
 \bigl(\mathcal D_{H,1}^n+2\mathcal D_{H,2}^{n+1/2}
 +\mathcal D_{H,1}^{n+1}\bigr)+\frac{\dt}{4}\sum_{n=1}^{\ell-1}
 \bigl(\mathcal F_{H,1}^n+2\mathcal F_{H,2}^{n+1/2}
 +\mathcal F_{H,1}^{n+1}\bigr).
\end{aligned}
\label{eq:cross-helicity-cumulative}
\end{equation}
\end{remark}

\subsubsection{Els\"asser energies and Alfv\'enic alignment}
\label{subsubsec:discrete-alignment}
For $1\le n\le N$, define the time-centered total energy
\begin{equation}
	\mathcal E_h^n
	=\frac12E_2^n
	+\frac14\bigl(E_1^{n-1/2}+E_1^{n+1/2}\bigr).
	\label{eq:combined-energy}
\end{equation}
The two integer-time Els\"asser energies are defined by
\begin{equation}
	D_+^n=\mathcal E_h^n+\sqrt{s}\,H_c^n,
	\qquad
	D_-^n=\mathcal E_h^n-\sqrt{s}\,H_c^n.
	\label{eq:Ds-plus-minus}
\end{equation}
They admit the nonnegative sum-of-squares representation
\begin{equation}
	\begin{aligned}
		D_\pm^n=\frac18\bigl(
		\bigl\|\bm u_1^{n-1/2}\pm\sqrt{s}\bm B_1^n\bigr\|^2
		+\bigl\|\bm u_1^{n+1/2}\pm\sqrt{s}\bm B_1^n\bigr\|^2+\bigl\|\bm u_2^n\pm\sqrt{s}\bm B_2^{n-1/2}\bigr\|^2
		+\bigl\|\bm u_2^n\pm\sqrt{s}\bm B_2^{n+1/2}\bigr\|^2
		\bigr).
	\end{aligned}
	\label{eq:Dminus}
\end{equation}
 In particular, $D_\pm^n=0$ implies the following
Alfv\'enic alignment relations:
\begin{equation}
\begin{aligned}
 \bm u_1^{n-1/2}=\bm u_1^{n+1/2}
 =\overline{\bm u}_1^n=\mp\sqrt{s}\bm B_1^n,\quad\bm u_2^n=\mp\sqrt{s}\bm B_2^{n-1/2}
 =\mp\sqrt{s}\bm B_2^{n+1/2}
 =\mp\sqrt{s}\overline{\bm B}_2^n.
\end{aligned}
 \label{eq:Dminus-averaged-alignment}
\end{equation}

\begin{corollary}[Discrete Els\"asser energies and Alfv\'enic alignment in the ideal case]
	\label{prop:alfvenic-alignment}
If $Re=Rm=\infty$ and $\bm f=\mathbf0$, then, for either choice of sign,
\begin{equation}
 D_\pm^n=D_\pm^1,\qquad 1\le n\le N.
 \label{eq:discrete-elsasser-conservation}
\end{equation}
In particular, if $D_\pm^1=0$ for either choice of sign, then
the relations in \eqref{eq:Dminus-averaged-alignment}
hold with the corresponding signs for $1\le n\le N$.
\end{corollary}
\begin{proof}
In the ideal, unforced case, Theorem~\ref{thm:discrete-energy-cross}
gives conservation of both energy sequences and $H_c^n$. Thus
\eqref{eq:combined-energy} gives $\mathcal E_h^n=\mathcal E_h^1$,
and \eqref{eq:Ds-plus-minus} yields
\eqref{eq:discrete-elsasser-conservation}. The exact alignment
statements follow from the nonnegative sum of squares in
\eqref{eq:Dminus}.
\end{proof}

\begin{remark}[Els\"asser-energy balances with dissipation]
\label{rem:elsasser-dissipative-balances}
Write $\nu=Re^{-1}$ and $\eta=Rm^{-1}$. Combining
\eqref{eq:continuous-energy-laws} and \eqref{eq:continuous-cross-law},
the continuous Els\"asser energies in
\eqref{eq:continuous-elsasser-energies} satisfy
\begin{equation}
\begin{aligned}
 \frac{\mathrm dD_\pm}{\mathrm dt}
 &=-\frac{\nu}{2}\sum_{i=1}^2\|\boldsymbol\omega_i\|^2
 -\frac{s\eta}{2}\sum_{i=1}^2\|\bm J_i\|^2\mp\frac{\sqrt{s}(\nu+\eta)}{2}
 \sum_{i=1}^2(\boldsymbol\omega_i,\bm J_i)
 +\frac12\sum_{i=1}^2(\bm f,\bm u_i\pm\sqrt{s}\bm B_i).
\end{aligned}
 \label{eq:continuous-elsasser-balance}
\end{equation}
Similarly, \eqref{eq:energy1-law}, \eqref{eq:energy2-law}, and
\eqref{eq:cross-helicity-law} give, for $1\le n\le N-1$,
\begin{equation}
\begin{aligned}
 \frac{D_\pm^{n+1}-D_\pm^n}{\dt}
 &=-\frac{\nu}{4}\bigl(
 \sum_{k=n}^{n+1}\|\overline{\boldsymbol\omega}_2^k\|^2
 +2\|\overline{\boldsymbol\omega}_1^{n+1/2}\|^2\bigr)-\frac{s\eta}{4}\bigl(
 \sum_{k=n}^{n+1}\|\overline{\bm J}_1^k\|^2
 +2\|\overline{\bm J}_2^{n+1/2}\|^2\bigr)\\
 &\quad\mp\frac{\sqrt{s}\nu}{4}\bigl(
 \sum_{k=n}^{n+1}(\overline{\boldsymbol\omega}_2^k,\bm J_2^k)
 +2(\overline{\boldsymbol\omega}_1^{n+1/2},\bm J_1^{n+1/2})\bigr)\mp\frac{\sqrt{s}\eta}{4}\bigl(
 \sum_{k=n}^{n+1}(\boldsymbol\omega_1^k,\overline{\bm J}_1^k)
 +2(\boldsymbol\omega_2^{n+1/2},\overline{\bm J}_2^{n+1/2})\bigr)\\
 &\quad+\frac14\sum_{k=n}^{n+1}
 (\bm f^k,\overline{\bm u}_1^k\pm\sqrt{s}\bm B_1^k)
 +\frac12(\bm f^{n+1/2},\overline{\bm u}_2^{n+1/2}
 \pm\sqrt{s}\bm B_2^{n+1/2}).
\end{aligned}
 \label{eq:discrete-elsasser-balance}
\end{equation}
The terms in \eqref{eq:discrete-elsasser-balance} are grouped in the same
order as in \eqref{eq:continuous-elsasser-balance}: vorticity dissipation,
current-density dissipation, vorticity--current cross terms, and forcing.
The two cross terms retain the distinct staggered-time pairings
associated with viscosity and resistivity.
In the unforced case with $\nu=\eta=\kappa$, the continuous balance becomes
$\mathrm dD_\pm/\mathrm dt
=-\tfrac{\kappa}{2}\sum_{i=1}^2
\|\boldsymbol\omega_i\pm\sqrt{s}\bm J_i\|^2\le0$.
The discrete balance does not by itself assert monotonicity of $D_\pm^n$.

%Applying Jensen's inequality to \eqref{eq:Dminus} gives, for either choice of sign,
%\begin{equation}
%	\|\overline{\bm u}_1^n\pm\sqrt{s}\bm B_1^n\|^2
%	+\|\bm u_2^n\pm\sqrt{s}\overline{\bm B}_2^n\|^2
%	\le4D_\pm^n.
%	\label{eq:time-centered-alignment-defect-bound}
%\end{equation}
%The Els\"asser energies $D_\pm^n$ thus provide upper bounds for the
%corresponding time-centered field differences.
\end{remark}

The magnitudes of $D_\pm^n$ control the absolute differences between
the corresponding velocity and scaled magnetic fields in the $L^2$ norm,
but they scale quadratically with the field amplitudes.
To assess relative alignment, we use
\eqref{eq:Ds-plus-minus} to define, for $\mathcal E_h^n>0$, the
normalized cross helicity
\begin{equation}
 \sigma_c^n=\frac{\sqrt{s}\,H_c^n}{\mathcal E_h^n}
 =\frac{\Dp^n-\Dm^n}{\Dp^n+\Dm^n},
 \qquad -1\leq\sigma_c^n\leq1.
 \label{eq:normalized-cross-helicity}
\end{equation}
Dividing the Els\"asser energies by the total energy gives the
dimensionless relative departures
\begin{equation}
 \frac{D_-^n}{\mathcal E_h^n}=1-\sigma_c^n,
 \qquad
 \frac{D_+^n}{\mathcal E_h^n}=1+\sigma_c^n.
 \label{eq:normalized-alignment-defects}
\end{equation}
$\sigma_c^n$ reflects both directional correlation and relative field
amplitudes, with values near $1$ and $-1$ indicating near global
Alfv\'enic alignment and anti-alignment, respectively. Values near zero indicate weak net
correlation, which may result from cancellation
between locally aligned and anti-aligned regions. Conserving both
$\mathcal E_h^n$ and $H_c^n$ preserves $\sigma_c^n$, whereas conserving
energy alone does not.

\subsubsection{Magnetic- and fluid-helicity balances}
\label{subsubsec:magnetic-helicity}
\label{subsubsec:fluid-helicity}

\noindent\textit{Magnetic potentials and helicity.}
To reconstruct the magnetic potentials, define
\begin{equation*}
	Z_{E,h}=\left\{\bm e_h\in\Eh:
	\left(\bm e_h,\nabla q_h\right)=0\quad\forall q_h\in\Vh\right\},
	\qquad
	Z_{F,h}=\left\{\bm f_h\in\Fh:
	\diver\bm f_h=0\right\}.
\end{equation*}
Exactness of \eqref{eq:discrete-exact-complex} gives an isomorphism
\begin{equation*}
	\curl:Z_{E,h}\to Z_{F,h}.
\end{equation*}
By the adjoint relation \eqref{eq:weakcurl}, $\curlh:Z_{F,h}\to Z_{E,h}$ is
the corresponding adjoint isomorphism.

These isomorphisms ensure that the gauge-fixed potentials are uniquely determined by
\begin{alignat}{2}
	\bm A_1^{n+1/2}&\in Z_{E,h},
	&\qquad \curl\bm A_1^{n+1/2}&=\bm B_2^{n+1/2},
	\label{eq:potential-A1-curl}\\
	\bm A_2^n&\in Z_{F,h},
	&\qquad \curlh\bm A_2^n&=\bm B_1^n.
	\label{eq:potential-A2-curl}
\end{alignat}
To compute magnetic helicity, we reconstruct the potentials from the magnetic fields by constrained linear solves.

For $1\le n\le N$, define two integer-time magnetic helicities
\begin{align*}
	H_{M,1}^n=\bigl(\overline{\bm{A}}_1^n,\bm{B}_1^n\bigr), \qquad
	H_{M,2}^n=\bigl(\bm{A}_2^n,\overline{\bm{B}}_2^n\bigr).
\end{align*}
Using \eqref{eq:weakcurl} and
\eqref{eq:potential-A1-curl}--\eqref{eq:potential-A2-curl}, we obtain,
for either choice of sign,
\begin{equation}
	\bigl(\bm{A}_1^{n\pm1/2},\bm{B}_1^n\bigr)
	=\bigl(\curl \bm{A}_1^{n\pm1/2},\bm{A}_2^n\bigr)
	=\bigl(\bm{B}_2^{n\pm1/2},\bm{A}_2^n\bigr).
	\label{eq:Hm-adjacent-equality}
\end{equation}
Averaging the two sign choices in \eqref{eq:Hm-adjacent-equality} yields
\begin{equation*}
	H_{M,1}^n=\bigl(\curl\overline{\bm{A}}_1^n,\bm{A}_2^n\bigr)
	=\bigl(\overline{\bm{B}}_2^n,\bm{A}_2^n\bigr)=H_{M,2}^n.
\end{equation*}
We denote their common value by $H_M^n=H_{M,1}^n=H_{M,2}^n$.

\noindent\textit{Fluid helicity.}
For $1\le n\le N$, define at integer time
\begin{align*}
	H_{F,1}^n=\bigl(\overline{\bm{u}}_1^n,\boldsymbol{\omega}_1^n\bigr),\qquad
	H_{F,2}^n=\bigl(\bm{u}_2^n,\overline{\boldsymbol{\omega}}_2^n\bigr).
\end{align*}
The adjoint curl identity and the discrete vorticity relations give,
for either choice of sign,
\begin{equation}
 \bigl(\bm u_1^{n\pm1/2},\boldsymbol\omega_1^n\bigr)
 =\bigl(\curl\bm u_1^{n\pm1/2},\bm u_2^n\bigr)
 =\bigl(\boldsymbol\omega_2^{n\pm1/2},\bm u_2^n\bigr).
 \label{eq:Hf-adjacent-equality}
\end{equation}
Averaging the two sign choices in \eqref{eq:Hf-adjacent-equality} gives
$H_{F,1}^n=H_{F,2}^n$. We denote their common value by $H_F^n=H_{F,1}^n=H_{F,2}^n$.

\begin{theorem}[Discrete magnetic- and fluid-helicity balances]
\label{thm:discrete-helicity-balances}
\label{thm:discrete-magnetic-balance}
\label{thm:discrete-fluid-balance}
For $1\le n\le N-1$, the common helicities satisfy:

\textup{(i) Magnetic helicity.}
\begin{equation}
	H_M^{n+1}-H_M^n=
	-\frac{\dt}{2}\bigl(
	\mathcal D_{M,1}^n+2\mathcal D_{M,2}^{n+1/2}
	+\mathcal D_{M,1}^{n+1}\bigr),
	\label{eq:magnetic-helicity-one-step}
\end{equation}
Thus $Rm=\infty$ implies $H_M^{n+1}=H_M^n$.

\textup{(ii) Fluid helicity.}
\begin{equation}
H_F^{n+1}-H_F^n
 =-\frac{\dt}{2}\bigl(
 \mathcal D_{F,1}^n+2\mathcal D_{F,2}^{n+1/2}
 +\mathcal D_{F,1}^{n+1}\bigr)
 +\frac{\dt}{2}\bigl(
 \mathcal F_{F,1}^n+2\mathcal F_{F,2}^{n+1/2}
 +\mathcal F_{F,1}^{n+1}\bigr).
\label{eq:fluid-helicity-one-step}
\end{equation}
In the unforced inviscid hydrodynamic case,
$\bm B_1=\bm B_2=\mathbf0$, $Re=\infty$, and
$\bm f=\mathbf0$, one has $H_F^{n+1}=H_F^n$.\\
For $1\le n\le N$, the magnetic-helicity contributions are defined by
\begin{align*}
	\mathcal D_{M,1}^n&=Rm^{-1}(\overline{\bm J}_1^n,\bm B_1^n),\\
	\mathcal D_{M,2}^{n-1/2}&=Rm^{-1}(\overline{\bm J}_2^{n-1/2},\bm B_2^{n-1/2}),
\end{align*}
and the fluid-helicity contributions by
\begin{align*}
	\mathcal D_{F,1}^n
	&=Re^{-1}\bigl(\overline{\boldsymbol\omega}_2^n,
	\curl\boldsymbol\omega_1^n\bigr)
	-s\bigl(\overline{\bm J}_1^n\times\bm B_1^n,
	\boldsymbol\omega_1^n\bigr),\\
	\mathcal D_{F,2}^{n-1/2}
	&=Re^{-1}\bigl(\curl\overline{\boldsymbol\omega}_1^{n-1/2},
	\boldsymbol\omega_2^{n-1/2}\bigr)
	-s\bigl(\overline{\bm J}_2^{n-1/2}\times\bm B_2^{n-1/2},
	\boldsymbol\omega_2^{n-1/2}\bigr),\\
	\mathcal F_{F,1}^n&=\bigl(\bm f^n,\boldsymbol\omega_1^n\bigr), \qquad
	\mathcal F_{F,2}^{n-1/2}
	=\bigl(\bm f^{n-1/2},\boldsymbol\omega_2^{n-1/2}\bigr).
\end{align*}
The fluid-helicity terms include the Lorentz contribution.
\end{theorem}
\begin{proof}
\textup{(i)} \textit{Magnetic helicity.}
Test \eqref{eq:half-B} with $\bm C_2=\bm A_2^n$ and
\eqref{eq:int-B} with $\bm C_1=\bm A_1^{n+1/2}$.
The potential relations $\curlh\bm A_2^n=\bm B_1^n$ and
$\curl\bm A_1^{n+1/2}=\bm B_2^{n+1/2}$,
together with the scalar triple-product identity, eliminate the induction terms:
\begin{align*}
 (\overline{\bm u}_1^n\times\bm B_1^n,\curlh\bm A_2^n)
 &=(\overline{\bm u}_1^n\times\bm B_1^n,\bm B_1^n)=0,\\
 (\overline{\bm u}_2^{n+1/2}\times\bm B_2^{n+1/2},\curl\bm A_1^{n+1/2})
 &=(\overline{\bm u}_2^{n+1/2}\times\bm B_2^{n+1/2},\bm B_2^{n+1/2})=0.
\end{align*}
Thus, multiplying the tested equations by $\dt$ yields
\begin{align}
 (\bm A_2^n,\bm B_2^{n+1/2}-\bm B_2^{n-1/2})
 +\dt\mathcal D_{M,1}^n&=0,
 \label{eq:Hm-half-bridge}\\
 (\bm A_1^{n+1/2},\bm B_1^{n+1}-\bm B_1^n)
 +\dt\mathcal D_{M,2}^{n+1/2}&=0.
 \label{eq:Hm-int-bridge}
\end{align}
Average \eqref{eq:Hm-half-bridge} at $n$ and $n+1$ and add the result
to \eqref{eq:Hm-int-bridge}. By
\eqref{eq:Hm-adjacent-equality} and the definition of
$H_M^n$, the inner-product increments combine to $H_M^{n+1}-H_M^n$,
which proves \eqref{eq:magnetic-helicity-one-step}.
If $Rm=\infty$, both $\mathcal D_{M,i}$ vanish, so $H_M^{n+1}=H_M^n$.

\smallskip
\noindent\textup{(ii)} \textit{Fluid helicity.}
Test \eqref{eq:half-u} with $\bm v_1=\boldsymbol\omega_1^n$ and
\eqref{eq:int-u} with $\bm v_2=\boldsymbol\omega_2^{n+1/2}$.
The scalar triple-product identity eliminates the convection terms:
\[
 (\boldsymbol\omega_1^n\times\overline{\bm u}_1^n,\boldsymbol\omega_1^n)=0,
 \qquad
 (\boldsymbol\omega_2^{n+1/2}\times\overline{\bm u}_2^{n+1/2},
 \boldsymbol\omega_2^{n+1/2})=0.
\]
The curl relations $\boldsymbol\omega_1^n=\curlh\bm u_2^n$ and
$\boldsymbol\omega_2^{n+1/2}=\curl\bm u_1^{n+1/2}$,
together with $\curl\nabla=0$ and $\diver\curl=0$, eliminate the pressure terms:
\[
 (\nabla p_0^n,\boldsymbol\omega_1^n)
 =(\curl\nabla p_0^n,\bm u_2^n)=0,
 \qquad (p_3^{n+1/2},\diver\boldsymbol\omega_2^{n+1/2})=0.
\]
Multiplying the tested equations by $\dt$ therefore gives
\begin{align}
 (\bm u_1^{n+1/2}-\bm u_1^{n-1/2},\boldsymbol\omega_1^n)
 +\dt\mathcal D_{F,1}^n&=\dt\mathcal F_{F,1}^n,
 \label{eq:Hf-half-bridge}\\
 (\bm u_2^{n+1}-\bm u_2^n,\boldsymbol\omega_2^{n+1/2})
 +\dt\mathcal D_{F,2}^{n+1/2}&=\dt\mathcal F_{F,2}^{n+1/2}.
 \label{eq:Hf-int-bridge}
\end{align}
Similarly, average \eqref{eq:Hf-half-bridge} at $n$ and $n+1$ and add
the result to \eqref{eq:Hf-int-bridge}.
By \eqref{eq:Hf-adjacent-equality} and the definition of $H_F^n$,
the inner-product increments combine to
$H_F^{n+1}-H_F^n$, proving \eqref{eq:fluid-helicity-one-step}.
If $\bm B_1=\bm B_2=\bm0$, $Re=\infty$, and $\bm f=\bm0$,
all $\mathcal D_{F,i}$ and $\mathcal F_{F,i}$ vanish, so $H_F^{n+1}=H_F^n$.
\end{proof}

\begin{remark}[Cumulative magnetic- and fluid-helicity balances]
\label{rem:magnetic-fluid-cumulative}
Summing the one-step laws gives the cumulative balances, for
$2\le\ell\le N$:
\begin{equation}
	H_M^\ell-H_M^1=
	-\dfrac{\dt}{2}\sum_{n=1}^{\ell-1}
	\bigl(\mathcal D_{M,1}^n+2\mathcal D_{M,2}^{n+1/2}
	+\mathcal D_{M,1}^{n+1}\bigr).
	\label{eq:magnetic-helicity-cumulative}
\end{equation}
\begin{equation}
H_F^\ell-H_F^1
 +\frac{\dt}{2}\sum_{n=1}^{\ell-1}
 \bigl(\mathcal D_{F,1}^n+2\mathcal D_{F,2}^{n+1/2}
 +\mathcal D_{F,1}^{n+1}\bigr)
 -\frac{\dt}{2}\sum_{n=1}^{\ell-1}
 \bigl(\mathcal F_{F,1}^n+2\mathcal F_{F,2}^{n+1/2}
 +\mathcal F_{F,1}^{n+1}\bigr)=0.
\label{eq:fluid-helicity-cumulative}
\end{equation}
\end{remark}

\subsubsection{Energy lower bounds from magnetic helicity}
\label{subsubsec:energy-lower-bounds}
The discrete Poincar\'e estimates
\begin{equation}
	\norm{\bm A_1}\le C_1\norm{\bm B_2},
	\qquad
	\norm{\bm A_2}\le C_1\norm{\bm B_1},
	\label{eq:poincare}
\end{equation}
hold with a mesh-independent constant $C_1>0$ for a fixed polynomial order.
The first is the discrete Poincar\'e inequality on the gauge-fixed space
$Z_{E,h}$; see \cite{MR2269741,MR2594630,MR2009375,MR2317829} for the underlying
theory. The zero-mean reduction and orthogonality to the curl kernel remove
the harmonic and gradient components.
For the second estimate, exactness gives $\bm e_h\in Z_{E,h}$ such that
$\curl\bm e_h=\bm A_2\in Z_{F,h}$. The same discrete Poincar\'e inequality
implies $\|\bm e_h\|\le C_1\|\bm A_2\|$.
Using the adjoint relation \eqref{eq:weakcurl} and
$\curlh\bm A_2=\bm B_1$, we obtain
\[
\begin{aligned}
 \|\bm A_2\|^2
 &=\bigl(\bm A_2,\curl\bm e_h\bigr)
 =\bigl(\curlh\bm A_2,\bm e_h\bigr)
 =\bigl(\bm B_1,\bm e_h\bigr)\le\|\bm B_1\|\,\|\bm e_h\|
 \le C_1\|\bm B_1\|\,\|\bm A_2\|.
\end{aligned}
\]
This proves the second estimate with the same constant $C_1$.
These estimates yield the energy lower bound below,
a discrete analogue of the Arnold inequality \cite{MR1612569}.

\begin{corollary}[Energy lower bounds]
\label{thm:energy-lower-bounds}
If $Rm=\infty$, then for $1\le n\le N$,
\begin{equation}
 E_2^nE_1^{n+1/2}\ge\frac{s^2}{4C_1^2}|H_M^n|^2.
 \label{eq:productlower}
\end{equation}
If also $\bm f=\mathbf0$ and the conserved value $H_M=H_M^1$ is nonzero, both staggered energies have positive lower bounds:
\begin{equation}
 E_1^{n+1/2}\ge C_1^*|H_M|^2>0,\qquad
 E_2^n\ge C_2^*|H_M|^2>0.
 \label{eq:strictlower}
\end{equation}
Here $0<Re\le\infty$,
$K_M=s^2/(4C_1^2)$, $C_1^*=K_M/E_2^0$, and
$C_2^*=K_M/E_1^{1/2}$.
Neither energy can decay to zero.
\end{corollary}
\begin{proof}
Since $Rm=\infty$, \eqref{eq:Hm-half-bridge} gives
$(\bm A_2^n,\bm B_2^{n+1/2}-\bm B_2^{n-1/2})=0$.
Hence the definition $H_M^n=(\bm A_2^n,\overline{\bm B}_2^n)$ reduces to
$H_M^n=(\bm A_2^n,\bm B_2^{n+1/2})$.
The Cauchy--Schwarz inequality and \eqref{eq:poincare} yield
\begin{equation}
 |H_M^n|\le\norm{\bm A_2^n}\bigl\|\bm B_2^{n+1/2}\bigr\|
 \le C_1\norm{\bm B_1^n}\bigl\|\bm B_2^{n+1/2}\bigr\|.
 \label{eq:Hm-Poincare}
\end{equation}
The magnetic parts of \eqref{eq:E1}--\eqref{eq:E2} give
$\norm{\bm B_1^n}\le\sqrt{2E_2^n/s}$ and
$\bigl\|\bm B_2^{n+1/2}\bigr\|\le\sqrt{2E_1^{n+1/2}/s}$.
Substitution into \eqref{eq:Hm-Poincare} gives
$|H_M^n|\le(2C_1/s)\sqrt{E_2^nE_1^{n+1/2}}$, proving \eqref{eq:productlower}.

If also $\bm f=\bm0$ and $H_M\ne0$, the energy laws imply
$E_1^{n+1/2}\le E_1^{1/2}$ and $E_2^n\le E_2^0$.
Magnetic-helicity conservation gives $H_M^n=H_M$, so \eqref{eq:productlower} becomes
\begin{equation}
 E_2^nE_1^{n+1/2}\ge K_M|H_M|^2>0.
 \label{eq:positive-product}
\end{equation}
Using the upper bounds in \eqref{eq:positive-product} gives
$E_1^{n+1/2}\ge (K_M/E_2^0)|H_M|^2=C_1^*|H_M|^2$ and
$E_2^n\ge (K_M/E_1^{1/2})|H_M|^2=C_2^*|H_M|^2$, respectively.
Both constants depend only on the initial data, $s$, and $C_1$, and are independent
of $n$. This proves \eqref{eq:strictlower}.
\end{proof}

\section{Numerical experiments}
\label{sec:numerical-experiments}
Section~\ref{subsec:convergence-test} assesses spatial accuracy,
Section~\ref{sec:cross-role} examines the role of cross helicity,
Section~\ref{subsec:conservation-tests} checks the conservation properties,
and Section~\ref{sec:orszag-tang} illustrates planar vortex evolution.
Table~\ref{tab:results-guide} identifies the conservation laws tested below.

\subsection{Convergence tests}
\label{subsec:convergence-test}

\begin{figure}[!htbp]
	\centering
	\begin{subfigure}[t]{0.31\textwidth}
		\includegraphics[width=\linewidth]{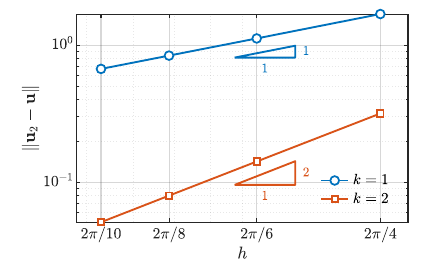}
		\caption{$\|\bm u_2^N-\bm u^N\|$.}
	\end{subfigure}\hfill
	\begin{subfigure}[t]{0.31\textwidth}
		\includegraphics[width=\linewidth]{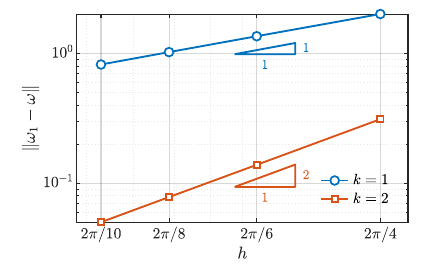}
		\caption{$\|\boldsymbol\omega_1^N-\boldsymbol\omega^N\|$.}
	\end{subfigure}\hfill
	\begin{subfigure}[t]{0.31\textwidth}
		\includegraphics[width=\linewidth]{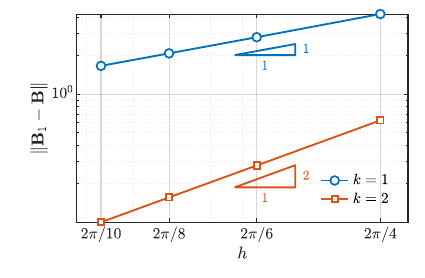}
		\caption{$\|\bm B_1^N-\bm B^N\|$.}
	\end{subfigure}
	\makebox[\textwidth][c]{%
	\begin{subfigure}[t]{0.31\textwidth}
		\includegraphics[width=\linewidth]{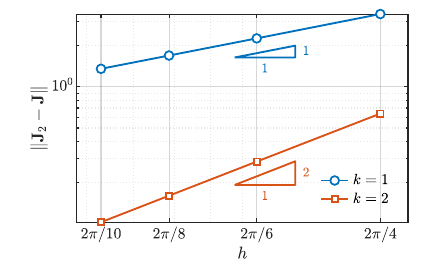}
		\caption{$\|\bm J_2^N-\bm J^N\|$.}
	\end{subfigure}\hspace{0.06\textwidth}%
	\begin{subfigure}[t]{0.31\textwidth}
		\includegraphics[width=\linewidth]{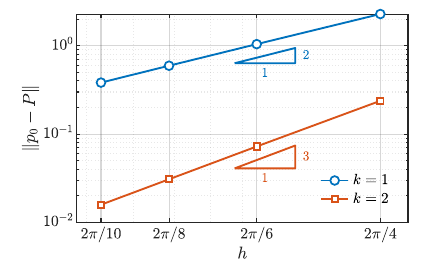}
		\caption{$\|p_0^N-P^N\|$.}
	\end{subfigure}}
	\caption{Final-time errors on the integer grid.}
	\label{fig:convergence-final-integer}
\end{figure}

\begin{figure}[!htbp]
	\centering
	\begin{subfigure}[t]{0.31\textwidth}
		\includegraphics[width=\linewidth]{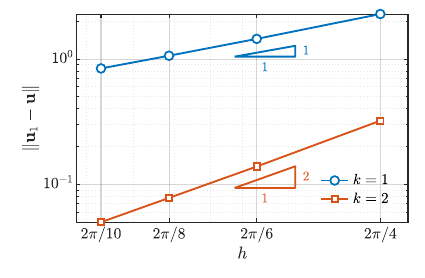}
		\caption{$\|\bm u_1^{N+1/2}-\bm u^{N+1/2}\|$.}
	\end{subfigure}\hfill
	\begin{subfigure}[t]{0.31\textwidth}
		\includegraphics[width=\linewidth]{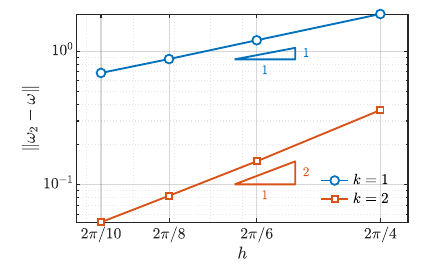}
		\caption{$\|\boldsymbol\omega_2^{N+1/2}-\boldsymbol\omega^{N+1/2}\|$.}
	\end{subfigure}\hfill
	\begin{subfigure}[t]{0.31\textwidth}
		\includegraphics[width=\linewidth]{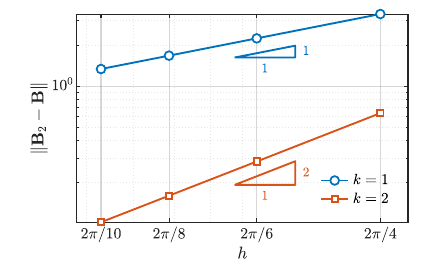}
		\caption{$\|\bm B_2^{N+1/2}-\bm B^{N+1/2}\|$.}
	\end{subfigure}
	\makebox[\textwidth][c]{%
	\begin{subfigure}[t]{0.31\textwidth}
		\includegraphics[width=\linewidth]{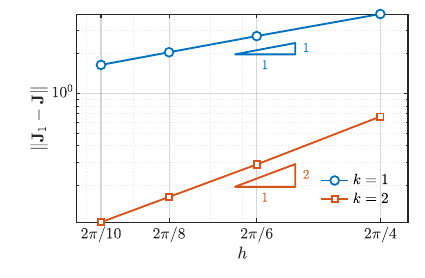}
		\caption{$\|\bm J_1^{N+1/2}-\bm J^{N+1/2}\|$.}
	\end{subfigure}\hspace{0.06\textwidth}%
	\begin{subfigure}[t]{0.31\textwidth}
		\includegraphics[width=\linewidth]{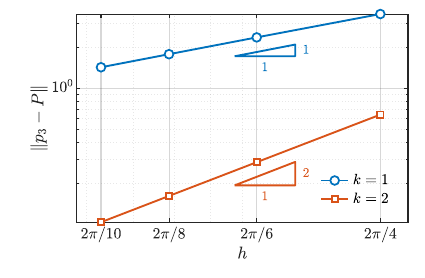}
		\caption{$\|p_3^{N-1/2}-P^{N-1/2}\|$.}
	\end{subfigure}}
	\caption{Last available errors on the half-integer grid.}
	\label{fig:convergence-final-half}
\end{figure}

On $\Omega=(0,2\pi)^3$, use
\[
\bm{u}=\mathrm{e}^{-t}
\begin{pmatrix}
	\sin y\\
	0\\
	\dfrac12\cos y
\end{pmatrix},
\qquad
\bm{B}=\mathrm{e}^{-t}
\begin{pmatrix}
	\sin y\\
	0\\
	2\cos y
\end{pmatrix},
\qquad
P=\sin x,
\qquad
\bm{f}=
\begin{pmatrix}
	\cos x\\
	-\dfrac{15}{4}\mathrm{e}^{-2t}\sin y\cos y\\
	0
\end{pmatrix}.
\]
We set $s=1$, $Re=Rm=1$, $T=1$, and $\dt=0.01$. The mesh sizes are $h=2\pi/m$ for $m=4,6,8,10$, and each cube is divided into six tetrahedra. The integer-grid fields are
\[
\bigl(
\bm{u}_2^n,
\boldsymbol{\omega}_1^n,
\bm{B}_1^n,
\bm{J}_2^n
\bigr),
\]
and the half-integer-grid fields are
\[
\bigl(
\bm{u}_1^{n+1/2},
\boldsymbol{\omega}_2^{n+1/2},
\bm{B}_2^{n+1/2},
\bm{J}_1^{n+1/2}
\bigr).
\]
Here $N=T/\dt=100$. The errors in these integer- and half-integer-grid
fields are evaluated at $t_N=1$ and $t_{N+1/2}=1.005$, respectively.
The pressures are $p_0^n$ on the integer grid and $p_3^{n+1/2}$ on the
half-integer grid. Their errors use the last available values:
$p_0^N$ at $t_N=1$ and $p_3^{N-1/2}$ at $t_{N-1/2}=0.995$, respectively.
We take $k=1$ and $k=2$. Figures~\ref{fig:convergence-final-integer} and~\ref{fig:convergence-final-half} show the errors at these time levels. The reference triangles indicate $O(h^k)$ for the vector fields and $p_3$, and $O(h^{k+1})$ for $p_0$.

\subsection{Role of cross helicity}
\label{sec:cross-role}

We compare the proposed and Mao--Xi methods \cite{MR4917731} to investigate the role of the cross-helicity balance in Alfv\'enic alignment. An ideal test examines the preservation of an initially aligned state, while a weakly dissipative test examines alignment dynamics from non-aligned initial data.

Both tests use $\Omega=(0,1)^3$, with $k=1$ and six tetrahedra per Cartesian cube.
In each test, the two methods use the same mesh $h=1/12$, time-step size, physical parameters, and initial data.
In the figure captions, $N=12$ denotes mesh subdivisions per coordinate direction, not time steps. Both tests are unforced.

For the Mao--Xi method, the centered cross-helicity diagnostic is defined by
\begin{equation}
 H_{c,\mathrm M}^n=\frac12\left[
 (\overline{\bm u}_{1,\mathrm M}^n,\overline{\bm B}_{1,\mathrm M}^n)
 +(\bm u_{2,\mathrm M}^n,\bm B_{2,\mathrm M}^n)\right].
 \label{eq:Mao-candidate-Hc}
\end{equation}
This quantity is not among the conservation laws established in \cite[Sec.~3.2]{MR4917731}.

\subsubsection{Ideal aligned test}
\label{subsubsec:aligned-ideal-test}

Set $\bm u_0=\bm B_0=U$, where
\[
 U=\begin{pmatrix}
 \sin2\pi z+\cos2\pi y\\
 \sin2\pi x+\cos2\pi z\\
 \sin2\pi y+\cos2\pi x
 \end{pmatrix},\qquad \curl U=2\pi U.
\]
We take $\dt=0.01$, $T=2$, $s=1$, and $Re=Rm=\infty$. By Corollary~\ref{thm:continuous-elsasser}, the initial relation
$\bm u=\bm B$ is preserved in this unforced ideal test
with $s=1$. The induction equation then gives $\partial_t\bm B=\mathbf0$, hence $\bm u(t)=\bm B(t)=U$.

\begin{figure}[htbp]
  \centering
  \begin{subfigure}[t]{0.32\textwidth}
    \includegraphics[width=\linewidth]{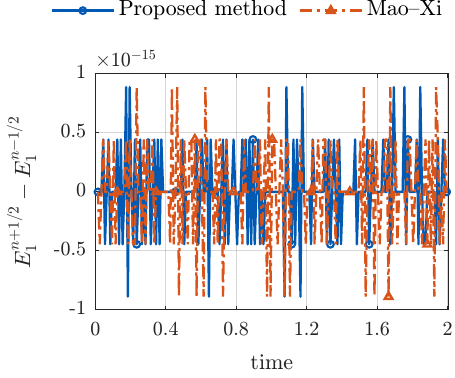}
    \caption{Increment of $E_1$ (half-integer steps).}
    \label{fig:aligned-E1}
  \end{subfigure}\hfill
  \begin{subfigure}[t]{0.32\textwidth}
    \includegraphics[width=\linewidth]{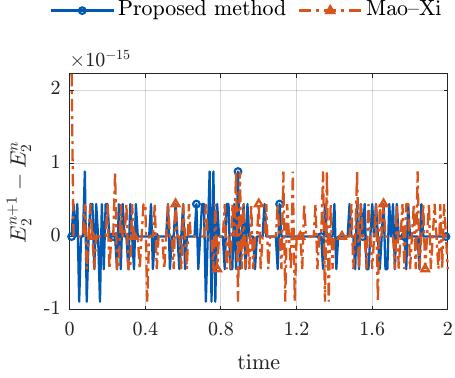}
    \caption{Increment of $E_2$ (integer steps).}
    \label{fig:aligned-E2}
  \end{subfigure}\hfill
  \begin{subfigure}[t]{0.32\textwidth}
    \includegraphics[width=\linewidth]{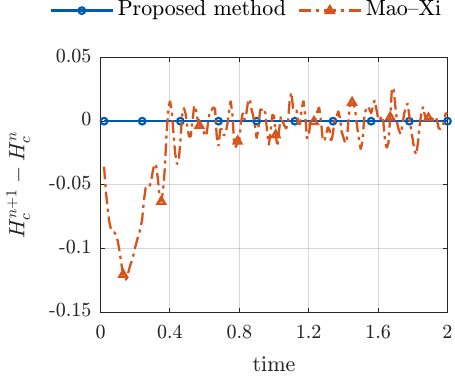}
    \caption{Increment of $H_c$.}
    \label{fig:aligned-Hc}
  \end{subfigure}
  \caption{Ideal aligned test: increments of $E_1$, $E_2$, and $H_c$ at $N=12$.}
  \label{fig:aligned-invariants}
\end{figure}

\begin{figure}[htbp]
  \centering
  \begin{subfigure}[t]{0.32\textwidth}
    \includegraphics[width=\linewidth]{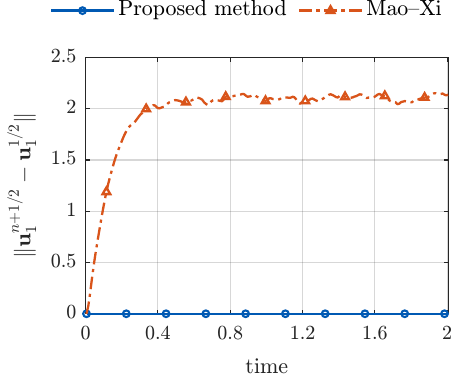}
    \caption{$\|\bm u_1^{n+1/2}-\bm u_1^{1/2}\|$.}
  \end{subfigure}\hspace{0.04\textwidth}%
  \begin{subfigure}[t]{0.32\textwidth}
    \includegraphics[width=\linewidth]{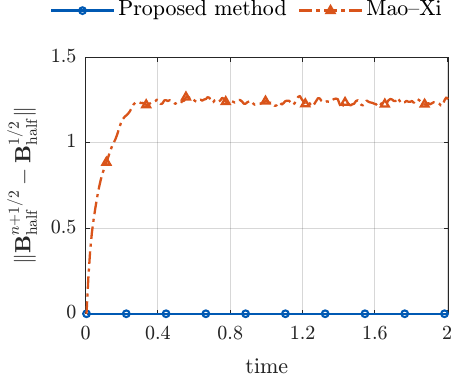}
    \caption{$\|\bm B_{\mathrm{half}}^{n+1/2}-\bm B_{\mathrm{half}}^{1/2}\|$.}
  \end{subfigure}

  \begin{subfigure}[t]{0.32\textwidth}
    \includegraphics[width=\linewidth]{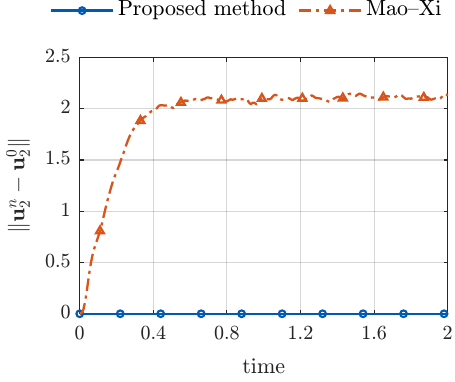}
    \caption{$\|\bm u_2^n-\bm u_2^0\|$.}
  \end{subfigure}\hspace{0.04\textwidth}%
  \begin{subfigure}[t]{0.32\textwidth}
    \includegraphics[width=\linewidth]{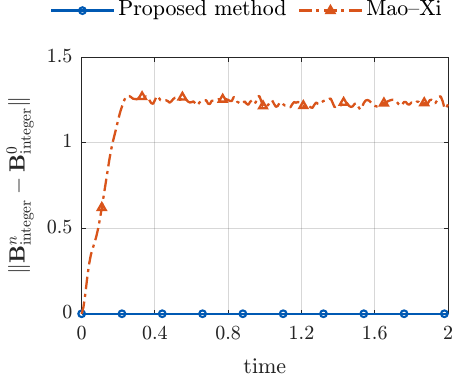}
    \caption{$\|\bm B_{\mathrm{integer}}^n-\bm B_{\mathrm{integer}}^0\|$.}
  \end{subfigure}
  \caption{Ideal aligned test: $L^2$ deviations from the initial discrete fields at $N=12$.
  For the proposed and Mao--Xi methods, respectively,
  $\bm B_{\mathrm{half}}$ denotes $\bm B_2$ and $\bm B_{1,\mathrm M}$,
  while $\bm B_{\mathrm{integer}}$ denotes $\bm B_1$ and $\bm B_{2,\mathrm M}$.}
  \label{fig:aligned-L2}
\end{figure}

Figure~\ref{fig:aligned-invariants} shows that both methods conserve their two energy sequences, but only the proposed method conserves cross helicity. Figure~\ref{fig:aligned-L2} compares each discrete field with its initialized value: $t_{1/2}$ for half-integer-step fields and $t_0$ for integer-step fields. The proposed method retains its initial discrete state, consistently with Corollary~\ref{prop:alfvenic-alignment} and the stationary continuous solution $\bm u(t)=\bm B(t)=U$ derived above. The Mao--Xi method departs from its initial state despite energy conservation.

\subsubsection{Weakly dissipative non-aligned test}
\label{subsubsec:nonaligned-dissipative-test}

This test examines the dissipative cross-helicity balance and velocity--magnetic alignment. Let $k=2\pi$ and set
\[
\bm u_0=\frac1{\sqrt3}
\begin{pmatrix}
	\sin(kz)+\cos(2ky)\\
	\sin(kx)+\cos(2kz)\\
	\sin(ky)+\cos(2kx)
\end{pmatrix},
\quad
\bm q_0=\frac1{\sqrt3}
\begin{pmatrix}
	\cos(kz)+\sin(2ky)\\
	\cos(kx)+\sin(2kz)\\
	\cos(ky)+\sin(2kx)
\end{pmatrix},
\quad
\bm B_0=\rho\bm u_0+\sqrt{1-\rho^2}\,\bm q_0,
\quad
\rho=0.6.
\]
We use $\dt=0.001$, $T=5$, $s=1$, and $Re^{-1}=Rm^{-1}=10^{-3}$.

\begin{figure}[!htbp]
  \centering
  \begin{subfigure}[t]{0.32\textwidth}
    \includegraphics[width=\linewidth]{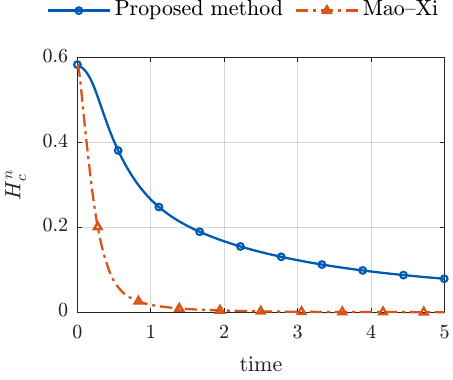}
    \caption{Cross helicity $H_c$.}
  \end{subfigure}\hspace{0.04\textwidth}%
  \begin{subfigure}[t]{0.32\textwidth}
    \includegraphics[width=\linewidth]{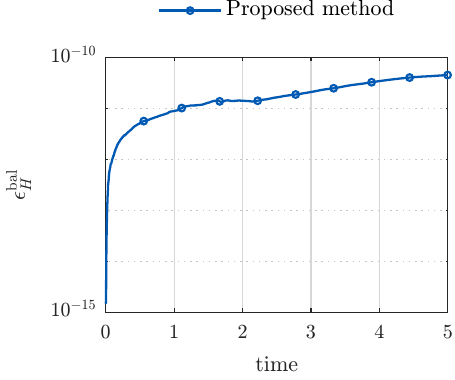}
    \caption{$\epsilon_H^{\mathrm{bal}}(t)$.}
  \end{subfigure}
  \caption{Weakly dissipative test: cross helicity and the cumulative balance residual at $N=12$.}
  \label{fig:dissipative-balance}
\end{figure}

Figure~\ref{fig:dissipative-balance}(a) shows distinct numerical evolutions:
the cross helicity rapidly decays towards zero for the Mao--Xi method,
whereas it remains appreciably positive for the proposed method at $T=5$.
To further verify that the proposed method satisfies the cumulative
cross-helicity balance \eqref{eq:cross-helicity-cumulative} in this
unforced test, define the normalized residual
\begin{equation}
\begin{aligned}
 \epsilon_H^{\mathrm{bal}}(t_n)=\frac{|H_c^n-H_c^1+\mathcal B_H^n|}{H_{\mathrm{ref}}}, \qquad
 H_{\mathrm{ref}}=\max\{|H_c^1|,\mathcal E_h^1\}, \qquad
 \mathcal B_H^n=\frac{\dt}{4}\sum_{m=1}^{n-1}
 \left(\mathcal D_{H,1}^{m}+2\mathcal D_{H,2}^{m+1/2}
 +\mathcal D_{H,1}^{m+1}\right).
\end{aligned}
\label{eq:dissipative-cross-residual}
\end{equation}
The terms $\mathcal D_{H,i}$ are defined in \eqref{eq:cross-dissipation-one}--\eqref{eq:cross-dissipation-two}.
Figure~\ref{fig:dissipative-balance}(b) shows that the residual remains
below $10^{-10}$, verifying the predicted balance to high precision.
No residual is plotted for the Mao--Xi method because the discrete
balance law established here does not apply to it.

\begin{figure}[!htbp]
  \centering
  \begin{subfigure}[t]{0.32\textwidth}
    \includegraphics[width=\linewidth]{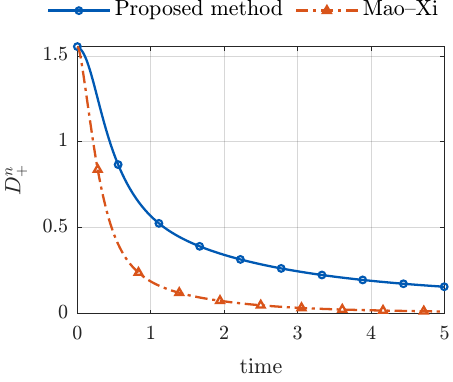}
    \caption{Els\"asser energy $D_+^n$.}
  \end{subfigure}\hfill
  \begin{subfigure}[t]{0.32\textwidth}
    \includegraphics[width=\linewidth]{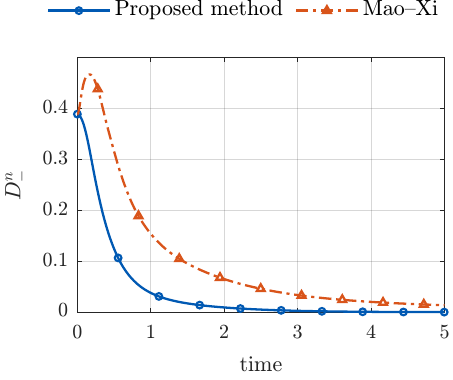}
    \caption{Els\"asser energy $D_-^n$.}
  \end{subfigure}\hfill
  \begin{subfigure}[t]{0.32\textwidth}
    \includegraphics[width=\linewidth]{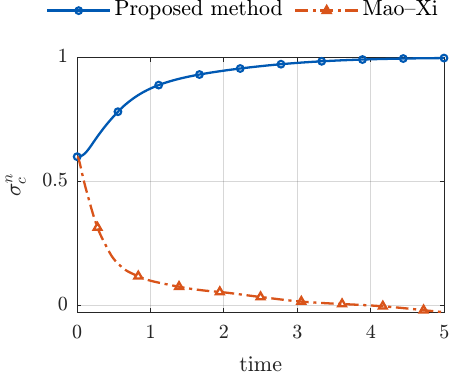}
    \caption{Normalized cross helicity $\sigma_c^n$.}
  \end{subfigure}
  \caption{Weakly dissipative test: $D_+^n$, $D_-^n$, and $\sigma_c^n$ at $N=12$.}
  \label{fig:dissipative-elsasser}
\end{figure}

Since cross helicity measures velocity--magnetic correlation, we next
examine whether the two methods produce different global alignment trends.
For the proposed method, the global diagnostics $D_{\pm}^n$ and
$\sigma_c^n$ are defined by \eqref{eq:Ds-plus-minus} and
\eqref{eq:normalized-cross-helicity}.
For the Mao--Xi method, use its two staggered
energies \cite[Theorem~3.5]{MR4917731} and define
\begin{align*}
 \mathcal E_{\mathrm M}^n=\frac14(E_{1,\mathrm M}^{n-1/2}+E_{1,\mathrm M}^{n+1/2})
 +\frac12E_{2,\mathrm M}^n, \qquad
 D_{\pm,\mathrm M}^n=\mathcal E_{\mathrm M}^n\pm\sqrt{s}\,H_{c,\mathrm M}^n, \qquad
 \sigma_{c,\mathrm M}^n=\frac{\sqrt{s}\,H_{c,\mathrm M}^n}{\mathcal E_{\mathrm M}^n}.
\end{align*}
Here $E_{1,\mathrm M}$ is the half-integer-step energy associated with
$(\bm u_{1,\mathrm M},\bm B_{1,\mathrm M})$, and
$E_{2,\mathrm M}$ is the integer-step energy associated with
$(\bm u_{2,\mathrm M},\bm B_{2,\mathrm M})$.
In the plots, each method uses its own energies and cross helicity
to compute $D_\pm^n$ and $\sigma_c^n$.
Figure~\ref{fig:dissipative-elsasser}(a,b) compares the two Els\"asser
energies. For the proposed method, $D_-^n$ becomes much smaller than
$D_+^n$, whereas the two remain comparable for the Mao--Xi method.
Accordingly, panel~(c) shows that $\sigma_c^n$ increases towards one
for the proposed method, indicating a tendency towards positive
Alfv\'enic alignment, while it decreases to values near zero for the
Mao--Xi method.

\begin{figure}[!htbp]
  \centering
  \begin{subfigure}[t]{0.32\textwidth}
    \includegraphics[width=\linewidth]{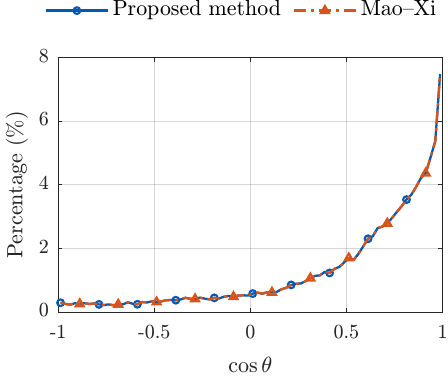}
    \caption{$t=0$.}
  \end{subfigure}\hspace{0.04\textwidth}%
  \begin{subfigure}[t]{0.32\textwidth}
    \includegraphics[width=\linewidth]{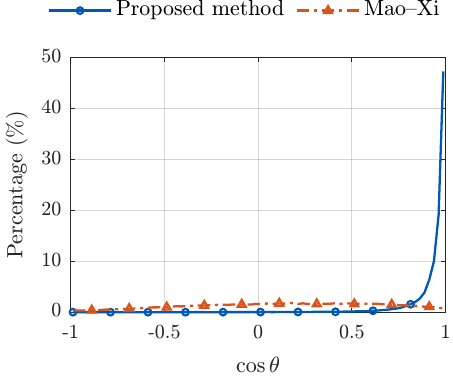}
    \caption{$t=1.25$.}
  \end{subfigure}

  \begin{subfigure}[t]{0.32\textwidth}
    \includegraphics[width=\linewidth]{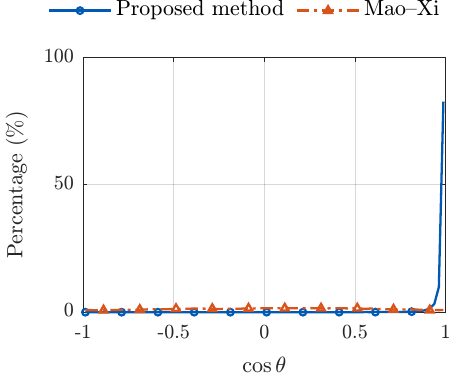}
    \caption{$t=3$.}
  \end{subfigure}\hspace{0.04\textwidth}%
  \begin{subfigure}[t]{0.32\textwidth}
    \includegraphics[width=\linewidth]{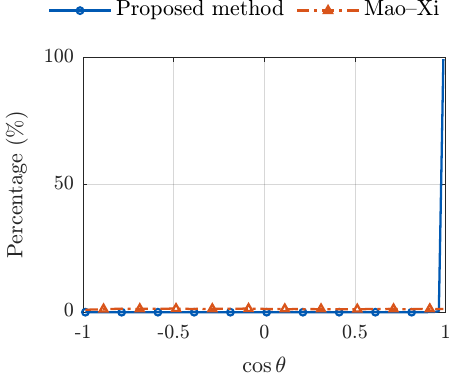}
    \caption{$t=5$.}
  \end{subfigure}
  \caption{Weakly dissipative test: distribution of the local alignment at $N=12$.}
  \label{fig:dissipative-alignment}
\end{figure}

To examine how these global trends are reflected locally, reconstruct the fields as
$\bm u_h^n=(\overline{\bm u}_1^n+\bm u_2^n)/2$ and
$\bm B_h^n=(\bm B_1^n+\overline{\bm B}_2^n)/2$, and define the angle
between them by
\begin{equation}
 \cos\theta(\bm x,t_n)=
 \frac{\bm u_h^n(\bm x)\cdot\bm B_h^n(\bm x)}
 {|\bm u_h^n(\bm x)|\,|\bm B_h^n(\bm x)|}.
 \label{eq:local-alignment-angle}
\end{equation}
For the Mao--Xi method, use
$\bm B_{h,\mathrm M}^n=(\overline{\bm B}_{1,\mathrm M}^n+\bm B_{2,\mathrm M}^n)/2$
with the same velocity reconstruction. At $t=0$, average the two projected initial approximations of each field.
For both methods, use the same $48^3$ uniformly spaced sampling points
$\bm x_\ell$ at $t=0,1.25,3,5$.
Each method retains points with
$|\bm u_h^n|\,|\bm B_h^n|>\tau_n$, where
\mbox{$\tau_n=\max\bigl\{10^{-14},\,10^{-12}\max_\ell
|\bm u_h^n(\bm x_\ell)|\,|\bm B_h^n(\bm x_\ell)|\bigr\}$}.
This excludes nearly vanishing denominators, for which the computed angle
can be unreliable. Divide $[-1,1]$ into 80 equal bins and plot the
percentage of retained points in each bin. Figure~\ref{fig:dissipative-alignment} shows that the proposed method's
distribution becomes concentrated near $\cos\theta=1$, indicating
increasing local alignment, while the Mao--Xi distribution remains
broader. This agrees with the different global alignment trends in
Figure~\ref{fig:dissipative-elsasser}.

The two methods produce distinct numerical solutions: the proposed
solution tends towards Alfv\'enic alignment, whereas the Mao--Xi solution
does not. This trend is qualitatively consistent with the two-dimensional
results of Ting et al.~\cite{ting1986}. Together with the discrete energy
and cross-helicity balances, this consistency lends further credibility
to the proposed numerical results.

\subsection{Conservation tests}
\label{subsec:conservation-tests}

All three-dimensional tests in this subsection use the lowest-order choice $k=1$ in the complex \eqref{eq:discrete-exact-complex}. The computational
domains are partitioned uniformly into cubes, and each cube is divided consistently into six tetrahedra to obtain a periodic tetrahedral mesh.

\subsubsection{Energy and cross-helicity conservation and divergence-free constraints}
\label{subsubsec:energy-cross-divergence-test}

\begin{figure}[htbp]
 \centering
 \begin{subfigure}[t]{0.32\textwidth}
  \includegraphics[width=\linewidth]{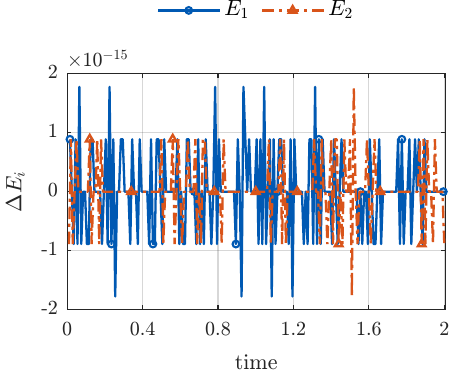}
  \caption{Energy increments $\Delta E_i$.}
 \end{subfigure}\hfill
 \begin{subfigure}[t]{0.32\textwidth}
  \includegraphics[width=\linewidth]{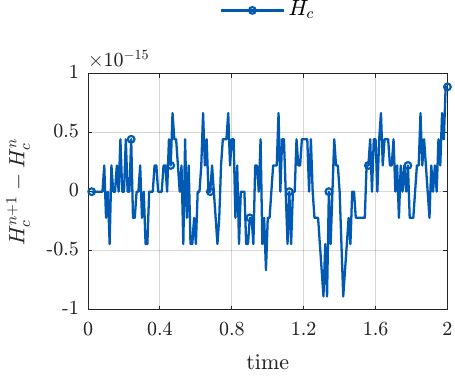}
  \caption{Cross-helicity increment.}
 \end{subfigure}\hfill
 \begin{subfigure}[t]{0.32\textwidth}
  \includegraphics[width=\linewidth]{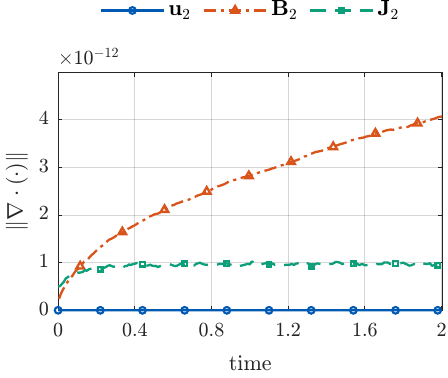}
  \caption{Divergence $L^2$ norms.}
 \end{subfigure}
 \caption{Energy and cross-helicity conservation and strong divergence-free constraints in the unforced ideal test.}
 \label{fig:ideal}
\end{figure}

This test examines energy and cross-helicity conservation in Theorem~\ref{thm:discrete-energy-cross}, together with the strong divergence-free constraints in Proposition~\ref{thm:discrete-divergence-preservation}. We use
$\Omega=(0,1)^3$, $s=1$, $\bm f=\mathbf0$, and $Re=Rm=\infty$, with the initial fields

\[
\bm{u}_0=
\begin{pmatrix}
	\sin 2\pi y+\cos 2\pi z+0.3\sin 4\pi y\\
	\sin 2\pi z+\cos 2\pi x+0.2\cos 4\pi z\\
	\sin 2\pi x+\cos 2\pi y+0.4\sin 4\pi x
\end{pmatrix},
\qquad
\bm{B}_0=
\begin{pmatrix}
	\sin 4\pi y+2\cos 2\pi z\\
	2\sin 2\pi z+\cos 4\pi x\\
	\sin 4\pi x-\cos 2\pi y
\end{pmatrix}.
\]
Both fields are periodic, divergence-free, and zero-mean. We take $h=1/8$, $\dt=0.01$, and $T=2$. The plotted increments are $\Delta E_1^n=E_1^{n+1/2}-E_1^{n-1/2}$, $\Delta E_2^n=E_2^{n+1}-E_2^n$, and $\Delta H_c^n=H_c^{n+1}-H_c^n$. Figure~\ref{fig:ideal}(a,b) shows that these increments remain at roundoff level. Panel~(c) plots the $L^2$ norms of $\diver\bm u_2$, $\diver\bm B_2$, and $\diver\bm J_2$; their near-zero values verify the strong divergence-free constraints.

\subsubsection{Magnetic-helicity conservation and energy lower bound}
\label{subsubsec:magnetic-helicity-test}

We conduct short- and long-time tests on $\Omega=(0,1)^3$, with $s=1$, $\bm f=\mathbf0$, $Re=1$, $Rm=\infty$, and the same initial fields as in Subsection~\ref{subsubsec:energy-cross-divergence-test}. The initial magnetic helicity is nonzero. The two equivalent expressions $H_{M,1}^n=H_{M,2}^n=H_M^n$ represent the same discrete magnetic helicity; the helicity panels plot their successive increments $H_{M,i}^{n+1}-H_{M,i}^n$, $i=1,2$.

\begin{figure}[!htbp]
 \centering
 \begin{subfigure}[t]{0.32\textwidth}
  \includegraphics[width=\linewidth]{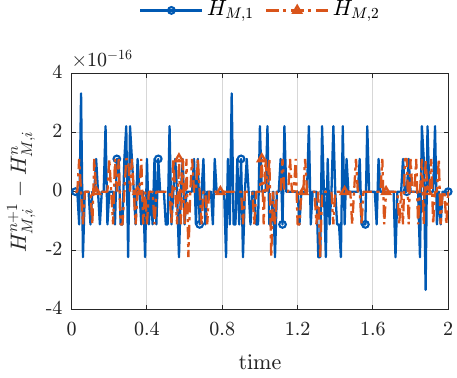}
  \caption{Magnetic-helicity increments.}
 \end{subfigure}\hfill
 \begin{subfigure}[t]{0.32\textwidth}
  \includegraphics[width=\linewidth]{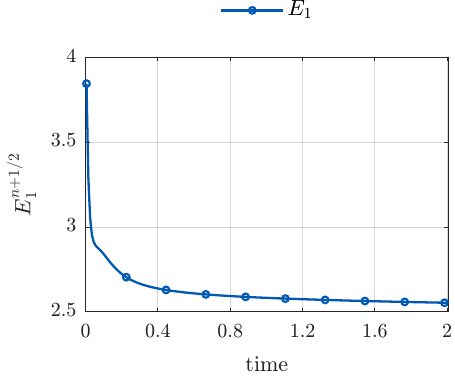}
  \caption{Half-integer-step energy $E_1$.}
 \end{subfigure}\hfill
 \begin{subfigure}[t]{0.32\textwidth}
  \includegraphics[width=\linewidth]{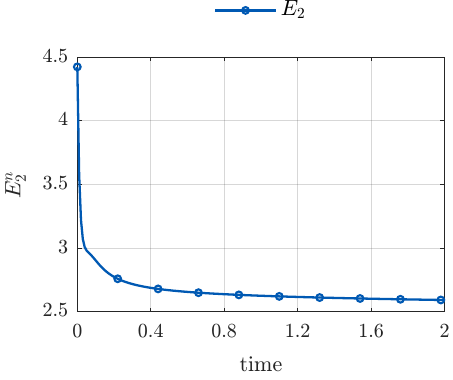}
  \caption{Integer-step energy $E_2$.}
 \end{subfigure}
 \par\smallskip
 \begin{subfigure}[t]{0.32\textwidth}
  \includegraphics[width=\linewidth]{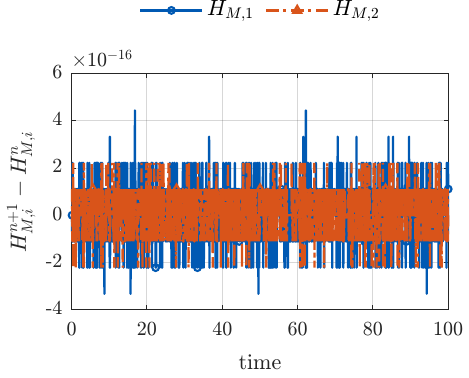}
  \caption{Magnetic-helicity increments.}
 \end{subfigure}\hfill
 \begin{subfigure}[t]{0.32\textwidth}
  \includegraphics[width=\linewidth]{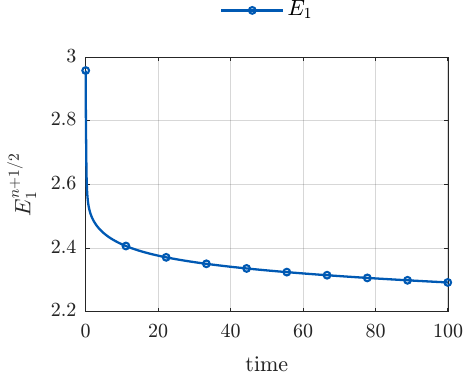}
  \caption{Half-integer-step energy $E_1$.}
 \end{subfigure}\hfill
 \begin{subfigure}[t]{0.32\textwidth}
  \includegraphics[width=\linewidth]{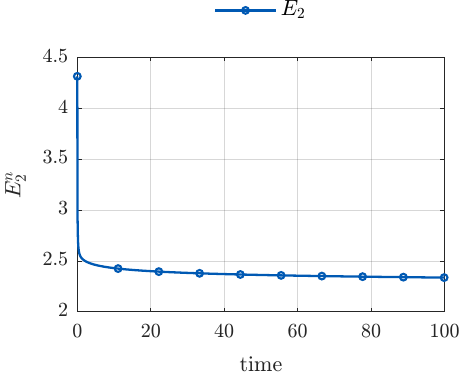}
  \caption{Integer-step energy $E_2$.}
 \end{subfigure}
 \caption{Magnetic-helicity conservation and energy evolution with $Re=1$ and
 $Rm=\infty$: (a--c) the short-time test with $h=1/8$, $\dt=0.01$, and $T=2$;
 (d--f) the long-time test with $h=1/6$, $\dt=0.05$, and $T=100$.}
 \label{fig:magnetic-comparison}
\end{figure}

For the short-time test, we take $h=1/8$, $\dt=0.01$, and $T=2$. Figure~\ref{fig:magnetic-comparison}(a) shows near-zero magnetic-helicity increments, while panels~(b,c) show the decrease of both energy sequences due to viscosity. Since $Rm=\infty$, magnetic helicity remains conserved despite this viscous energy dissipation.

The long-time test examines the energy evolution under the positive lower-bound constraint. We take $h=1/6$, $\dt=0.05$, and $T=100$. Figure~\ref{fig:magnetic-comparison}(d) again shows near-zero magnetic-helicity increments. The energy sequences in panels~(e,f) decrease but remain well above zero over the simulated time interval, consistently with the theoretical lower bounds. 
These tests illustrate magnetic-helicity conservation in Theorem~\ref{thm:discrete-helicity-balances}(i) and energy evolution consistent with the positive lower bounds in Corollary~\ref{thm:energy-lower-bounds}.

\subsubsection{Fluid-helicity conservation}
\label{subsubsec:fluid-helicity-test}

\begin{figure}[!htbp]
 \centering
 \includegraphics[width=0.32\textwidth]{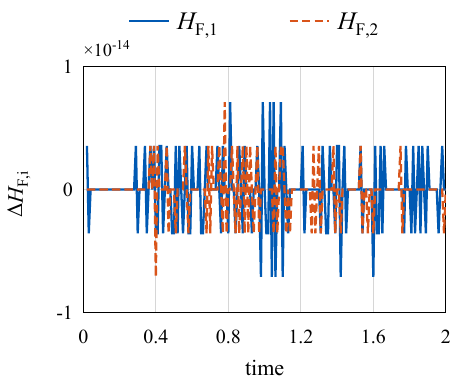}
 \caption{Increments of $H_{F,1}$ and $H_{F,2}$ in the inviscid hydrodynamic test.}
 \label{fig:fluid}
\end{figure}

This test examines fluid-helicity conservation in the hydrodynamic limit of Theorem~\ref{thm:discrete-helicity-balances}(ii). We use $\Omega=(0,1)^3$, $s=1$, and the same initial velocity $\bm u_0$ as in Subsection~\ref{subsubsec:energy-cross-divergence-test}, while setting
$\bm B_0=\mathbf0$, $\bm f=\mathbf0$, $Re=\infty$, and $Rm=1$.
We take $h=1/8$, $\dt=0.01$, and $T=2$. The magnetic fields remain zero, so the problem reduces to unforced inviscid hydrodynamics. The two equivalent diagnostics satisfy $H_{F,1}^n=H_{F,2}^n=H_F^n$ by \eqref{eq:Hf-adjacent-equality}. Figure~\ref{fig:fluid} plots their successive increments $H_{F,i}^{n+1}-H_{F,i}^n$, $i=1,2$; the near-zero values confirm fluid-helicity conservation.

\subsection{Orszag--Tang vortex}
\label{sec:orszag-tang}

The planar Orszag--Tang tests examine vortex evolution and current-sheet formation \cite{orszag1979}. Taking all fields to be independent of $z$ and setting the
$z$-components of velocity and magnetic field to zero reduces the three-dimensional equations exactly to two-dimensional incompressible MHD. The scalar curls are
\[
 \omega=\partial_xu_y-\partial_yu_x,\qquad
 J=\partial_xB_y-\partial_yB_x,
\]
and the divergence-free constraints become
$\partial_xu_x+\partial_yu_y=\partial_xB_x+\partial_yB_y=0$.
Both tests use $\Omega=(0,2\pi)^2$, periodic boundary conditions on all sides, and $\bm f=\mathbf0$. For these zero-mean planar magnetic fields, an out-of-plane vector potential gives $\bm A\cdot\bm B=0$, so the usual three-dimensional magnetic helicity is trivial. These tests therefore assess planar field evolution rather than a nontrivial magnetic-helicity invariant.

\subsubsection{Periodic Cartesian spaces}
\label{subsubsec:ot-compatible-complex}

Let $Q_{r,s}(\mathcal T_h)$ denote cellwise tensor-product polynomials of degrees at most $r$ and $s$ on a uniform periodic Cartesian mesh. Using $V,E,F,Q$ for the corresponding two-dimensional Sobolev spaces, set
\begin{alignat*}{2}
	V_h&=Q_{1,1}(\mathcal T_h)\cap V,
	\qquad E_h=\bigl[Q_{0,1}(\mathcal T_h)\times
	Q_{1,0}(\mathcal T_h)\bigr]\cap E,\\
	F_h&=\bigl[Q_{1,0}(\mathcal T_h)\times
	Q_{0,1}(\mathcal T_h)\bigr]\cap F,
	\qquad Q_h=Q_{0,0}(\mathcal T_h)\cap Q.
\end{alignat*}
Tangential and normal degrees of freedom are identified periodically. After removal of the constant modes, the two complexes are exact:
\begin{equation}
\begin{aligned}
 0&\longrightarrow V_h\xrightarrow{\gradop}E_h
       \xrightarrow{\curlop}Q_h\longrightarrow0,\\
 0&\longrightarrow V_h\xrightarrow{\gradop^\perp}F_h
       \xrightarrow{\divop}Q_h\longrightarrow0.
\end{aligned}
\label{eq:2d-discrete-exact-complexes}
\end{equation}
Here $\nabla^\perp q=(\partial_yq,-\partial_xq)$ and
$\curl(v_x,v_y)=\partial_xv_y-\partial_yv_x$. The first is the lowest-order tensor-product complex described in \cite[Sec.~3.2]{MR3670721}. The second is obtained by rotating the vector fields clockwise through $90^\circ$. The weak scalar curl is
\begin{equation}
 (\curlh\bm c_h,q_h)=(\bm c_h,\nabla^\perp q_h)
 \qquad\forall q_h\in V_h.
 \label{eq:2d-weak-curl}
\end{equation}
The field placement is
\begin{align*}
	&\bm u_1,\bm B_1\in E_h,\quad \omega_1,J_1,p_0\in V_h; \\
	&\bm u_2,\bm B_2\in F_h,\quad \omega_2,J_2,p_3\in Q_h.
\end{align*}

\subsubsection{\texorpdfstring{Viscous and resistive Orszag--Tang vortex}{Viscous and resistive Orszag--Tang vortex}}
\label{subsubsec:ot-dissipative}

\begin{figure}[!htbp]
 \centering
 \begin{subfigure}[t]{0.45\textwidth}
  \includegraphics[width=\linewidth]{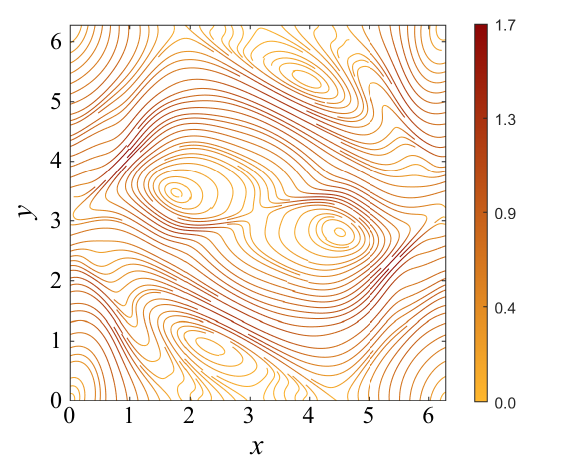}
  \caption{$\overline{\bm u}_1^n$.}
 \end{subfigure}\hfill
 \begin{subfigure}[t]{0.45\textwidth}
  \includegraphics[width=\linewidth]{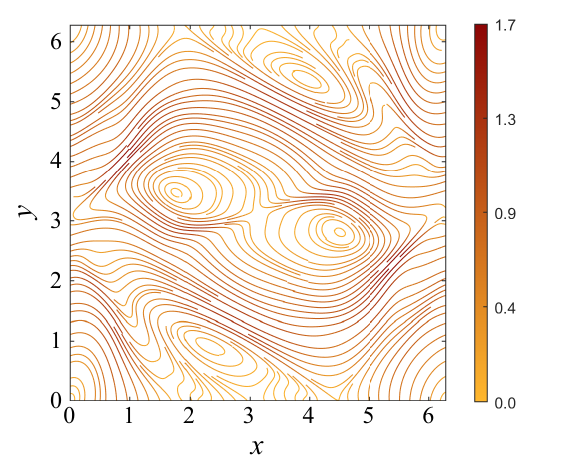}
  \caption{$\bm u_2^n$.}
 \end{subfigure}

 \begin{subfigure}[t]{0.45\textwidth}
  \includegraphics[width=\linewidth]{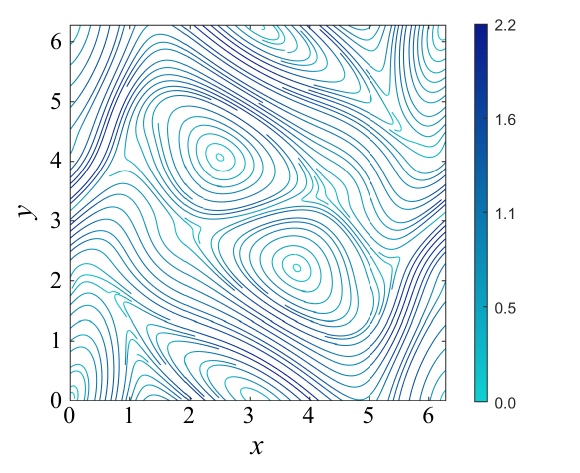}
  \caption{$\bm B_1^n$.}
 \end{subfigure}\hfill
 \begin{subfigure}[t]{0.45\textwidth}
  \includegraphics[width=\linewidth]{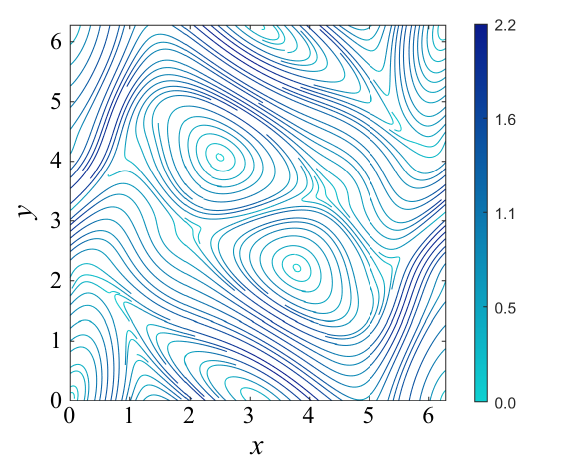}
  \caption{$\overline{\bm B}_2^n$.}
 \end{subfigure}
 \caption{Velocity streamlines (top row) and magnetic field lines (bottom row) for the viscous and resistive Orszag--Tang vortex at $t=2$. Colors indicate the magnitude of the corresponding field.}
 \label{fig:ot-2d-t2}
\end{figure}

\begin{figure}[!htbp]
 \centering
 \begin{subfigure}[t]{0.45\textwidth}
  \includegraphics[width=\linewidth]{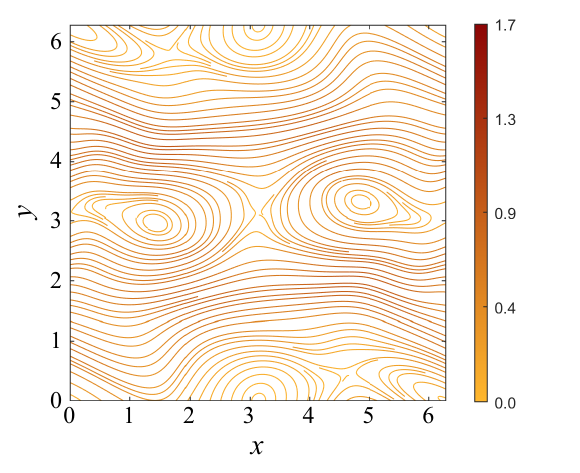}
  \caption{$\overline{\bm u}_1^n$.}
 \end{subfigure}\hfill
 \begin{subfigure}[t]{0.45\textwidth}
  \includegraphics[width=\linewidth]{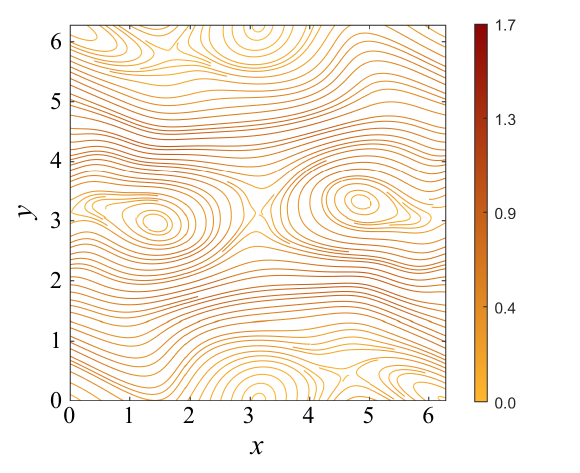}
  \caption{$\bm u_2^n$.}
 \end{subfigure}

 \begin{subfigure}[t]{0.45\textwidth}
  \includegraphics[width=\linewidth]{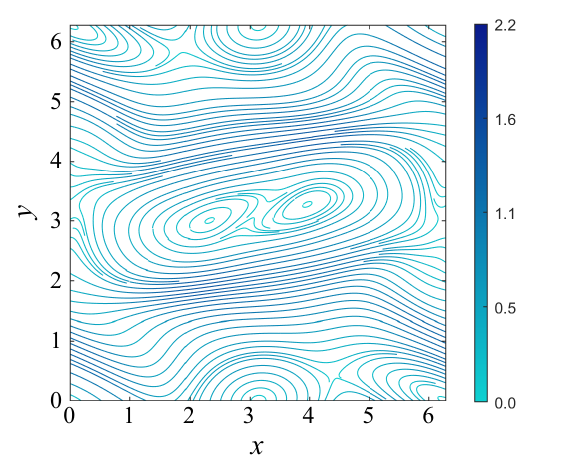}
  \caption{$\bm B_1^n$.}
 \end{subfigure}\hfill
 \begin{subfigure}[t]{0.45\textwidth}
  \includegraphics[width=\linewidth]{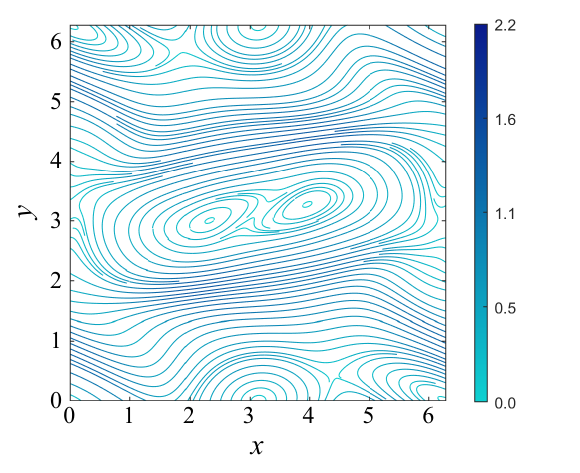}
  \caption{$\overline{\bm B}_2^n$.}
 \end{subfigure}
 \caption{Velocity streamlines (top row) and magnetic field lines (bottom row) for the viscous and resistive Orszag--Tang vortex at $t=10$. Colors indicate the magnitude of the corresponding field.}
 \label{fig:ot-2d-t10}
\end{figure}

The dissipative initial data are
\begin{equation*}
 \bm u^0=(-\sin y,\sin x),\qquad
 \bm B^0=(-\sin y,\sin 2x)
\end{equation*}
on $(0,2\pi)^2$.  We set $Re=Rm=100$, $s=1$,
$\bm f=\mathbf0$, $h=2\pi/200$, $\dt=0.01$, and $T=10$. 
The half-integer-time fields are averaged to integer times, so all four panels in each figure correspond to the same output time.
{The snapshots at $t=2$ and $t=10$ show the evolution and smoothing of the vortex structures.}
Compared with the fields at $t=2$ in Figure~\ref{fig:ot-2d-t2}, the velocity and magnetic fields at $t=10$ in Figure~\ref{fig:ot-2d-t10} vary more gradually in space. The two velocity approximations show similar spatial patterns, as do the two magnetic approximations. This evolution of the velocity and magnetic fields agrees qualitatively with the results in Section~5.5.1 of \cite{MR4917731} and the viscous--resistive results in Section~4.8 of \cite{MR4645615}.

\subsubsection{Nearly ideal Orszag--Tang vortex}
\label{subsubsec:ot-nearly-ideal}

Let
\begin{equation*}
 \phi^0=2\sin y-2\cos x,\qquad A^0=\cos(2y)-2\cos x,
\end{equation*}
and set $\bm u^0=\nabla^\perp\phi^0$ and
$\bm B^0=\nabla^\perp A^0$, namely
\begin{equation*}
 \bm u^0=(2\cos y,-2\sin x),\qquad
 \bm B^0=(-2\sin2y,-2\sin x).
\end{equation*}
We use $h=2\pi/256$, $\dt=0.001$, $T=0.7$, $s=1$, $\bm f=\mathbf0$, and $Re^{-1}=Rm^{-1}=10^{-6}$ as a small
regularization of the ideal and inviscid limit; the computed case therefore has small but nonzero viscosity and resistivity.

\begin{figure}[htbp]
 \centering
 \begin{subfigure}[t]{0.32\textwidth}
  \includegraphics[width=\linewidth]{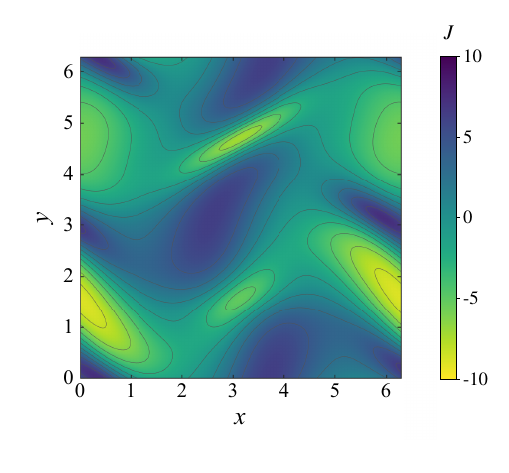}
  \caption{$\overline J_1^n$, $t=0.30$.}
 \end{subfigure}\hfill
 \begin{subfigure}[t]{0.32\textwidth}
  \includegraphics[width=\linewidth]{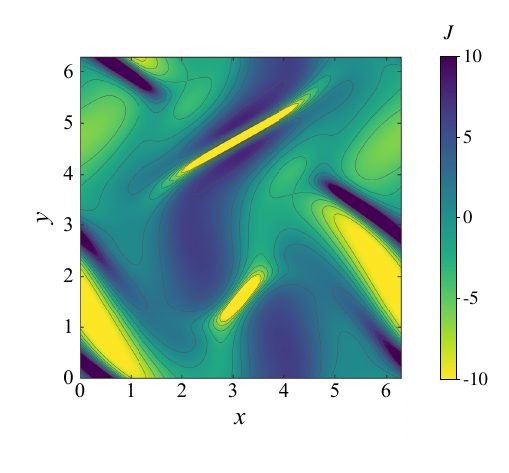}
  \caption{$\overline J_1^n$, $t=0.50$.}
 \end{subfigure}\hfill
 \begin{subfigure}[t]{0.32\textwidth}
  \includegraphics[width=\linewidth]{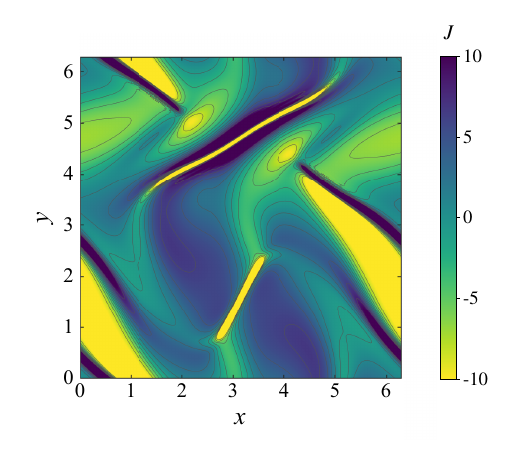}
  \caption{$\overline J_1^n$, $t=0.70$.}
 \end{subfigure}

 \begin{subfigure}[t]{0.32\textwidth}
  \includegraphics[width=\linewidth]{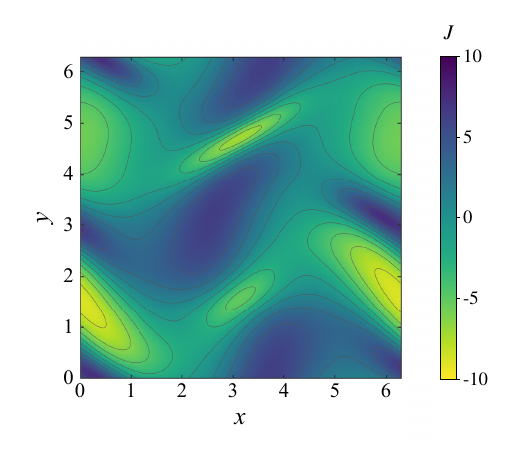}
  \caption{$J_2^n$, $t=0.30$.}
 \end{subfigure}\hfill
 \begin{subfigure}[t]{0.32\textwidth}
  \includegraphics[width=\linewidth]{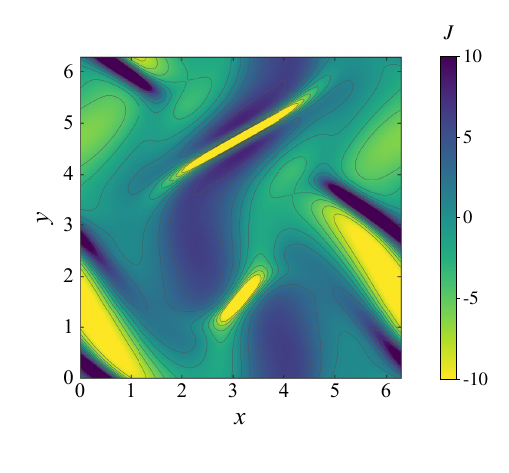}
  \caption{$J_2^n$, $t=0.50$.}
 \end{subfigure}\hfill
 \begin{subfigure}[t]{0.32\textwidth}
  \includegraphics[width=\linewidth]{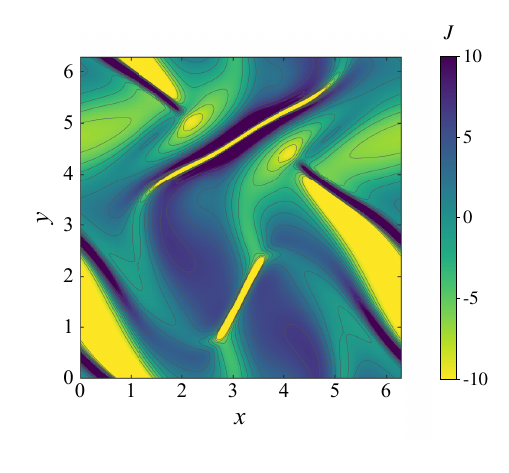}
  \caption{$J_2^n$, $t=0.70$.}
 \end{subfigure}
 \caption{Current densities $\overline J_1^n$ (top) and $J_2^n$ (bottom) for the nearly ideal Orszag--Tang vortex. All six panels use a common color scale ranging from $-10$ to $10$.}
 \label{fig:ot-ideal-current}
\end{figure}

At integer times, we display the two discrete approximations of the current density:
\begin{equation*}
	\overline J_1^n=\curlh\overline{\bm B}_2^n\in V_h,
	\qquad
	J_2^n=\curl\bm B_1^n\in Q_h.
\end{equation*}
Here $(\overline J_1^n,q_h)=(\overline{\bm B}_2^n,\nabla^\perp q_h)$
for all $q_h\in V_h$. Figure~\ref{fig:ot-ideal-current} shows the formation and progressive thinning of current sheets from $t=0.3$ to $T=0.7$. Both approximations show current sheets at similar locations. This evolution is qualitatively consistent with
the results reported in Section~5.3 of \cite{kraus2018} and Section~5.5.2 of \cite{MR4917731}.

\FloatBarrier

\FloatBarrier
\Needspace{13\baselineskip}
\section{Conclusions}
\label{sec:conclusions}

We have developed a dual-field finite element method for periodic incompressible MHD that reproduces the energy and cross-helicity balances through sequential linear staggered updates. The method also retains the strong and weak divergence-free constraints and the magnetic- and fluid-helicity balances.

In the unforced ideal case, simultaneous conservation of energy and cross helicity preserves both discrete Els\"asser energies and the normalized cross helicity. If either Els\"asser energy vanishes initially, its nonnegative sum-of-squares representation guarantees persistence of the corresponding exact Alfv\'enic relation. In the unforced non-resistive case, nonzero conserved magnetic helicity also gives positive lower bounds for both staggered energies.

Numerical tests show the expected spatial convergence rates and verify the discrete balance laws. The ideal aligned comparison with the Mao--Xi method demonstrates that energy conservation alone does not ensure preservation of the aligned state. For the weakly dissipative initial data considered, the proposed method satisfies the cumulative cross-helicity balance and develops increasing Alfv\'enic alignment, while the Mao--Xi method exhibits declining normalized cross helicity. Extending the formulation and its balance-law analysis to nonperiodic boundary conditions remains a direction for future work.

\section*{Declaration of competing interest}
The author declares no known competing financial interests or personal relationships that could have appeared to influence the work reported in this paper.

\section*{Data availability}
Data will be made available on request.

\begingroup
\footnotesize
\interlinepenalty=10000
\setlength{\bibsep}{1pt plus 0.2pt}
\bibliographystyle{mhd-references}
\bibliography{references}
\endgroup

\end{document}